\documentclass[twocolumn,english]{autart} 

\def\BibTeX{{\rm B\kern-.05em{\sc i\kern-.025em b}\kern-.08em
    T\kern-.1667em\lower.7ex\hbox{E}\kern-.125emX}}

\usepackage{epsfig}
\usepackage{subfigure}
\usepackage{amssymb,amsmath,amsfonts,layout,graphicx}
\usepackage{makeidx}
\usepackage{babel}
\usepackage{sublabel}
\usepackage{cases}
\usepackage{algorithm}

\usepackage{savesym}

\savesymbol{AND}
\savesymbol{OR}
\savesymbol{NOT}
\savesymbol{TO}
\savesymbol{COMMENT}
\savesymbol{BODY}
\savesymbol{IF}
\savesymbol{ELSE}
\savesymbol{ELSIF}
\savesymbol{FOR}
\savesymbol{WHILE}

\usepackage{algorithmic}

\usepackage
[
        letterpaper,
        left=1.91cm,
        right=1.91cm,
        top=1.91cm,
        bottom=2cm,
]
{geometry}

\newtheorem{lemma}{Lemma}
\newtheorem{assumption}{Assumption}
\newtheorem{proposition}{Proposition}
\newtheorem{theorem}{Theorem}
\newtheorem{definition}{Definition}
\newtheorem{corollary}{Corollary}
\newtheorem{remark}{Remark}

\newcommand{\R}{{\mathbb R}}

\DeclareMathOperator*{\argmin}{arg\,min}

\begin{document}
\begin{frontmatter}

\title{Connections between Mean-Field Game and Social Welfare Optimization}

\author[1]{Sen Li},   
\author[2]{Wei Zhang},
\author[2]{Lin Zhao}

\address[1]{Department of Mechanical Engineering, University of California, Berkeley, CA 94703, USA}  %
\address[2]{Department of Electrical and Computer Engineering, Ohio State University, Columbus, OH 43210, USA}

 \thanks{This work was supported in part by the National Science Foundation under grant CNS-1552838}
  \thanks{Email addresses: lisen1990@berkeley.edu (S. Li), zhang.491@osu.edu(W. Zhang), zhao.833@osu.edu(L. Zhao)}


\begin{abstract}
This paper studies the connection between a class of mean-field games and a social welfare optimization problem.  We consider a mean-field game in function spaces with a large population of agents, and each agent seeks to minimize an individual cost function. The cost functions of different agents are coupled through a mean-field term that depends on the mean of the population states.  
We show that although the mean-field game is not a potential game, under some mild condition the $\epsilon$-Nash equilibrium of the mean-field game coincides with the optimal solution to a social welfare optimization problem, and this is true even when the individual cost functions are non-convex. The connection enables us to evaluate and promote the efficiency of the mean-field equilibrium. In addition, it also leads to several important implications on the existence, uniqueness, and computation of the mean-field equilibrium. Numerical results are presented to validate the solution, and examples are provided to show the applicability of the proposed approach. 
\end{abstract}

\end{frontmatter}

\section{Introduction}
Mean-field games study the interactions among a large population of strategic agents, whose decision making is coupled through a mean-field term that depends on the statistical information of the overall population \cite{lasry2007mean}, \cite{lasry2006jeux}, \cite{lasry2006jeux1}, \cite{huang2007large}, \cite{gueant2011mean}. When the population size is large, each individual agent has a negligible impact on the mean-field term. This enables characterizing the game equilibrium via the interactions between the agent and the mean-field, instead of focusing on detailed interactions among all the agents. This idea was originally formalized in a series of seminal papers by Lasry and Lions  \cite{lasry2007mean}, \cite{lasry2006jeux} and by Huang et al.  \cite{huang2007large}, \cite{huang2006large}, where the mean-field equilibrium was characterized as the solution to an equation system that couples a backward Hamilton-Jacobi-Bellman equation and a forward Fokker-Planck-Kolmogorov equation. 
These seminal results attracted considerable research effort on various aspects of mean-field games. For instance, many works focused on  the existence \cite{huang2006large}, \cite{ferreira2015existence}, uniqueness \cite{lasry2007mean}, \cite{ahuja2016wellposedness}, and computation \cite{rammatico2016decentralized}, \cite{lachapelle2010computation}, \cite{achdou2010mean} of the mean-field equilibrium. For a more comprehensive review, please refer to \cite{lacker2015stochastic} and \cite{bensoussan2013mean}. Another strand of works extended the mean-field game model to more general settings such as heterogeneous agents \cite{huang2007large}, major-minor player model \cite{huang2010large}, \cite{bensoussan2015mean},  \cite{bensoussan2015mean2}, extended mean-field games \cite{gomes2013extended}, etc. Furthermore,  mean-field games also find abundant applications in economics \cite{gueant2009mean}, \cite{achdou2014partial}, crowd and population dynamics \cite{burger2013mean}, demand response \cite{couillet2012electrical}, \cite{bauso2014game}, \cite{ma2013decentralized}, networking \cite{huang2003individual}, coupled oscillators \cite{yin2012synchronization}, to name a few.

For many applications, it is important to analyze and quantify the efficiency of the mean-field equilibrium as compared to some socially optimal solutions. Along this line, the authors in \cite{huang2012social} and \cite{nourian2013nash} showed that the coordinator can design a mean-field game with an equilibrium that asymptotically achieves social optima as the population size goes to infinity. This result is true only when each agent in the game is {\em cooperative}. 
In the non-cooperative game setting, a recent work \cite{deori2017connection} showed that the Nash equilibrium of an electric vehicle charging game is socially optimal as the number of agents tends to infinity, under the assumption that the underlying game is a potential game. However, when the mean-field game is non-cooperative and does not admit a potential function, the mean-field equilibrium is shown to be {\em inefficient} in general.  
For instance, \cite{yin2014efficiency} employed a variational approach  to study the efficiency loss of mean-field equilibria for a synchronization game among oscillators. In \cite{balandat2013efficiency}, a mean-field congestion game was formulated, and numerical results were presented to show that the mean-field equilibrium is inefficient in general. 
In addition, \cite{nuno2015controlling} derived conditions under which the mean-filed equilibrium is efficient. However, since these conditions are quite restrictive. Instead of indicating the efficiency of the mean-field equilibrium, they are more useful in evaluating what happens in an inefficient equilibrium.

This paper studies the connection between mean-field games and social welfare optimization problems. We consider a class of mean-field games in vector spaces (potentially infinite-dimensional) with a large population of non-cooperative agents. Each agent seeks to minimize a cost functional. It couples the costs of other agents through a mean-field term that depends on the average of the population states. 
The key contribution of the paper lies in establishing the connection between the mean-field game and a {\em modified} social welfare optimization problem. This connection not only enables us to evaluate the promote the efficiency of the mean-field equilibrium, but also has several important implications on the existence, uniqueness, and computation of the mean-field equilibrium.  These contributions are summarized as follows: 
\begin{itemize}
\item First, we show that under some mild conditions, the mean-field equilibrium is  {\em efficient} with respect to  a {\em modified} social welfare optimization problem. This is true even when the individual costs and action spaces are non-convex. The connection enables us to {\em evaluate} and {\em promote} the efficiency of the mean-field equilibrium: given a mean-field game, we can construct the modified social welfare optimization problem to evaluate the efficiency of the equilibrium;  given a social cost function, we can design the mean-field game to provide a socially optimal mean-field equilibrium.

\item Second, we show that the mean-field equilibrium exists if the associated social welfare optimization problem has strong duality.  In addition, the other direction also holds under an additional monotonicity condition on the mean-field coupling term.  
Different from many existing works \cite{lasry2007mean}, \cite{ahuja2016wellposedness}, 
\cite{carmona2013probabilistic}, \cite{carmona2015probabilistic}, we obtain necessary and  sufficient conditions  for the existence of the mean-field equilibrium. This is stronger than the results in the existing literature.

\item Third, we show that the mean-field equilibrium is unique if the corresponding social welfare optimization problem is strictly convex. This recovers the results of some works on the uniqueness of the mean-field equilibrium \cite{lasry2007mean}, \cite{carmona2013probabilistic}, providing a novel interpretation of these results. 

\item Fourth, our result implies that computing the mean-field equilibrium is equivalent to solving a social welfare optimization problem. When this optimization problem is convex, various efficient algorithms can be employed to compute the mean-field equilibrium. According to our results, some existing methods on computing  
 the mean-field equilibrium  \cite{rammatico2016decentralized}, \cite{ma2013decentralized} can be interpreted as certain primal-dual algorithms in solving the associated social welfare optimization problem. To improve these algorithms, we provide an example of using the alternating direction method of multipliers \cite{boyd2011distributed} to compute the mean-field equilibrium. Simulation result shows that the proposed algorithm converges faster than the operator-based method.

\item Fifth, we discuss the differences between the modified social welfare  and the potential function in potential games. We show that in general, the mean-field game considered in this paper is not a potential game. In this case, the proposed method is the only way to connect the mean-field game to an optimization problem. Under additional assumptions(e.g., mean-field term is linear), the mean-field game  may reduce to a potential game. In this case, both the social welfare optimization and the potential function minimization are connected to the mean-field game. We show that compared to the potential function method,  the social welfare optimization 
(a)~provides a relaxed solution concept, (b) enables decentralized implementation and enjoys better scalability.
\end{itemize}

The rest of the paper proceeds as follows. The mean-field game is formulated in Section II. The  game equilibrium is characterized by a set of coupled equations in Section III. The connection between the mean-field equations and the social welfare optimization problem  is studied in Section IV. Section V discusses the implications on the existence, uniqueness and computation of the mean-field equilibrium. Section VI provides numerical results and examples.

\section{Mean-Field Games in Function Spaces}
This section formulates the class of mean-field games to be studied in this paper. Different from many existing works, we will describe the mean-field games in vector spaces (possibly infinite-dimensional). We show that our formulation includes many important classes of mean-field games as special cases \cite{huang2007large},  \cite{rammatico2016decentralized}, \cite{ma2013decentralized}, \cite{huang2010nce}. The vector space formulation allows us to more directly focus on aspects and challenges related to strategic interactions among the agents, instead of worrying about specific details and unnecessary technical conditions on system dynamics.

\subsection{The Mean-Field Game}
We consider a general mean-field game in a vector space among $N$ agents. Each agent $i$ is associated with a state variable $x_i$, a control input $u_i$ and a noise input $\pi_i$. The control input $u_i$ takes value in  $\mathcal{U}_i \subseteq \mathcal{U}$, where  $\mathcal{U}$ is an arbitrary vector space. The noise input $\pi_i$ is a random element in the measurable spaces $\Pi_i$ with the underlying probability space  $(\Omega, \mathcal{F}, P)$. The state of each agent is determined by the control and noise according to the following mapping $f_i:\mathcal{U}_i \times \Pi_i \rightarrow \mathcal{X}$:
\begin{equation}
\label{systemdynamics} 
x_i=f_i(u_i,\pi_i), \quad  u_i\in \mathcal{U}_i, 
\end{equation}
where $x_i$ is a random element that takes value in the space $\mathcal{X}$. When $x_i$ and $u_i$ are time-domain trajectories, the system dynamics are implicitly captured by (\ref{systemdynamics}).
To ensure that $x_i$ is well-defined, we impose the following assumption on $f_i(u_i,\pi_i)$:
\begin{assumption}
\label{assumptionf_i}
For each $u_i\in \mathcal{U}_i$, $f_i(u_i,\pi_i(\cdot)):\Omega \rightarrow \mathcal{X}$ is a measurable mapping with respect to $\mathcal{F}/\mathcal{Z}$, where $\mathcal{Z}$ is a $\sigma$-algebra on $\mathcal{X}$, 
\end{assumption}
Under Assumption \ref{assumptionf_i}, $x_i:\Omega \rightarrow \mathcal{X}$ is a measurable mapping with respect to $\mathcal{F}/\mathcal{Z}$. Therefore, $x_i$ is a well-defined random element that takes value in $\mathcal{X}$, with the expectation defined as the Bochner integral.

On the space $\mathcal{X}$, we define an inner product and a norm. The inner product is denoted as  $x\cdot y$ for $x,y\in \mathcal{X}$, and the norm is defined as $||x||=\sqrt{x\cdot x}$. We assume that $\mathcal{X}$ is complete.
\begin{assumption}
\label{assumptionf_i_2}
$\mathcal{X}$ is a Hilbert space. 
\end{assumption}
The completeness of $\mathcal{X}$ is mainly used to induce a special form of duality theory in vector spaces.  More elaborations can be found in Remark \ref{completeness}. 

Throughout the paper, we assume that $x_i$ and $x_j$ are independent. 
In addition, we assume the state $x_i$ has bounded second moment, i.e., there exists $C\geq 0$ such that $\mathbb{E} ||x_i||^2\leq C$ for all $i=1,\ldots,N$.  
In this case, the admissible control set can be defined as
$\bar{\mathcal{U}}_i=\{u_i\in \mathcal{U}_i \vert x_i=f_i(u_i,\pi_i), \mathbb{E} ||x_i||^2\leq C, x_i\in \mathcal{X}_i\}$, where $\mathcal{X}_i$ represents the local constraint on $x_i$, and $\mathcal{X}_i\subseteq \mathcal{X}$.

For each agent $i$, there is a cost function on the system state and control input. The costs of different agents are coupled through a mean-field term that depends on the average of the population state. We write it as follows:
\begin{equation}
\label{utilityfunction}
J_i(x_i,u_i,\bar{x})=\mathbb{E} \left( V_i(x_i,u_i)+F(\bar{x})\cdot  x_i+G(\bar{x})    \right),
\end{equation}
where $\bar{x}\in \mathcal{X}$ is the average of the population state, i.e., $\bar{x}\triangleq\dfrac{1}{N}\sum_{i=1}^N x_i$, $F:\mathcal{X}\rightarrow \mathcal{X}$  is the mean-field coupling term, and $G:\mathcal{X}\rightarrow\R$ is the cost associated with the mean-field term. The following conditions are imposed on $J_i$.
\begin{assumption}
\label{assumption_cost}
(i)  $F(\cdot)$  is globally Lipschitz continuous on $\mathcal{X}$ with constant $L$ , 
(ii)  $G(\cdot)$ is Fr\'echet differentiable on $\mathcal{X}$, and the gradient of $G(\cdot)$ at 0 is bounded, i.e., $||\nabla G(0)||<\infty$.
(iii) the gradient of $G(\cdot)$ is globally Lipschitz continuous on $\mathcal{X}$ with constant $\beta$, i.e., $||\nabla G(x)-\nabla G(y)||\leq \beta ||x-y||$, $\forall x,y\in \mathcal{X}$.  
\end{assumption}

The mean-field game considered in this paper is formulated as follows:
\begin{subequations}{\label{individualoptimization}}
\begin{align}
\min_{u_i} & \mathbb{E} \left(V_i(x_i,u_i)+F(\bar{x})\cdot  x_i +G(\bar{x})\right) \\ \label{constraintofmfg}
\text{s.t. } 
     &x_i=f_i(u_i, \pi_i), \quad u_i\in \bar{\mathcal{U}}_i,
\end{align}
\end{subequations}
where $\bar{x}\triangleq\dfrac{1}{N}\sum_{i=1}^N x_i$.  As  (\ref{individualoptimization})  is defined in vector spaces, it provides a unifying framework that incorporates a large body of literature. It includes both discrete-time \cite{rammatico2016decentralized}, \cite{ma2013decentralized} and continuous-time system \cite{huang2007large} as special cases, and addresses both deterministic and stochastic cases.  In the cost function (\ref{individualoptimization}), the inner product term $F(\bar{x})\cdot x_i$ can be either interpreted as the price multiplied by quantity \cite{lachapelle2010computation}, \cite{li2016on}, \cite{li2016social} or part of the quadratic penalty of the deviation of the system state from the population mean \cite{huang2007large}, \cite{huang2012social}. 
This structure is very common in the literature. We provide a few examples in the next subsection.

\begin{remark}
The proposed mean-field game has finite number of agents. This is different from many classic works \cite{lasry2007mean}, \cite{huang2006large}, where a continuum model is considered.  It is inspired by the continuum particle model in fluid dynamics. However, decision making problems with a continuum of decision-makers are uncommon in engineering applications.  In common cases, the continuum model is often used as an approximation of some finite system \cite{djehiche2016mean}. For such system, our finite model is directly applicable and naturally incorporates the heterogeneity of the agents. 
\end{remark}


\subsection{Examples}
\label{exampleofmeanfieldgame}
The mean-field game (\ref{individualoptimization}) captures a large class of important problems. In this subsection,  we present three examples:

\subsubsection{Discrete-time deterministic game}
First, consider a deterministic mean-field game in discrete-time \cite{rammatico2016decentralized}. The cost function of the game is  as follows:
\begin{equation}
\label{examplelqg1}
\min_{x_i} ||x_i||_Q^2+||x-\bar{x}||_{\Delta}^2+2(A\bar{x}+a)^T x_i
\end{equation}
where $x_i\in \mathbb{R}^T$ is the system state, $\bar{x}= \dfrac{1}{N} \sum_{i=1}^N x_i$ is the average state, $||x_i||_Q^2$ stands for $x_i^TQx_i$, $A\in \mathbb{R}^{T\times T}$, and $Q$ and $\Delta$ are symmetric positive definite matrices of appropriate
dimensions. Although  (\ref{examplelqg1}) is formulated in the  static form, it captures a class of finite-dimensional linear quadratic games, which can be transformed to (\ref{examplelqg1}) by plugging the linear dynamics in the cost function \cite{rammatico2016decentralized}.

Next, we show that the above game problem can be formulated as (\ref{individualoptimization}). Note that in (\ref{individualoptimization}),  the state space and the control space are both $\mathbb{R}^T$, and $x_i=u_i$, i.e., $f_i$ is the identity function. If we expand the norm and combine similar terms in (\ref{examplelqg1}), we can transform (\ref{examplelqg1}) to the following form:
\begin{equation}
\label{gameform1}
\min_{x_i} x_i^T(Q+\Delta)x_i+2a^Tx_i+2\bar{x}^T(A-\Delta)x_i+ \bar{x}^T \Delta \bar{x}
\end{equation}
Comparing (\ref{gameform1}) to (\ref{individualoptimization}), we have $V_i(x_i,u_i)=x_i^T(Q+\Delta)x_i+2a^Tx_i$, $F(\bar{x})=2(A-\Delta)\bar{x}$, and $G(\bar{x})=\bar{x}^T\Delta \bar{x}$. It is easy to verify that the game problem (\ref{gameform1}) satisfies Assumption \ref{assumptionf_i}-\ref{assumption_cost}. Therefore, (\ref{gameform1}) is a special case of the proposed mean-field game (\ref{individualoptimization}).

\subsubsection{Continuous-time stochastic game}
The second example is a linear quadratic Gaussian (LQG) game considered in Section II-A of \cite{huang2007large}:
\begin{align}
\label{examplelqg}
&\min_{\{u_i(t),t\geq 0\}} \mathbb{E} \int_0^\infty e^{-\rho t}\big[(x_i(t)-v(t))^2+ru_i(t)^2\big]dt \\
&\text{s.t.}
\begin{cases}
\label{linearsde}
dx_i(t)= u_i(t)dt+\sigma_i d\pi_i(t) \\
x_i(t)\in \mathbb{R},\quad u_i(t)\in \mathbb{R}, \quad i=1,\ldots,N,
\end{cases}
\end{align}
where $x_i(t)$ and $u_i(t)$ denote the state and control for the $i$th agent at time $t$, $\pi_i(t)$ is a standard scalar Brownian motion, $v(t)\triangleq b-c \dfrac{1}{N} \sum_{i=1}^N x_i(t)$ is the mean-field term, and $\rho, r,b, c>0$ are real constants. It is assumed that the Brownian motions $\pi_i=\{\pi_i(t),t\geq 0\}$ and $\pi_j=\{\pi_j(t),t\geq 0\}$ are independent for any $i\neq j$. 
Let $u_i=\{u_i(t),t\geq~0\}$. The authors in \cite{huang2007large} defined the admissible control set as $
\bar{\mathcal{U}}_i=\{u_i \vert  u_i(t)\text{ is adapted to the } \sigma-\text{algebra } \sigma(x_i(0),\pi_i(s),s\leq t), \text{ and } \mathbb{E} \int_{0}^{\infty} e^{-\rho t}[x_i(t)^2+u_i(t)^2]dt<\infty\}$. Let $x=\{x(t)\in \mathbb{R}, t\geq 0\}$, then the state space $\mathcal{X}$ is $\mathcal{X}=\{ x | \int_{[0,\infty)}e^{-\rho t}|x(t)|^2dt<\infty\}$. 

Next, we show that problem (\ref{examplelqg}) is a special case of the mean-field game (\ref{individualoptimization}). To this end, we  transform (\ref{examplelqg}) to the form of (\ref{individualoptimization}),  and verify Assumption \ref{assumptionf_i}- Assumption \ref{assumption_cost}. 

First,  the stochastic differential equation (\ref{linearsde}) can be solved explicitly. The solution is as follows:
\begin{equation}
\label{solutiontosde}
x_i(t)=\int_0^t u_i(s)ds+\sigma_i \int_0^t d\pi_i(s).
\end{equation}
It indicates that $x_i$ is determined by $\{u_i(s),0\leq s\leq t\}$, thus $f_i$ can be defined according to (\ref{solutiontosde}). Clearly, $x_i(t)$ is well-defined, and Assumption \ref{assumptionf_i} holds. 

Second, for any $x\in \mathcal{X}$ and $y \in \mathcal{X}$, we define their inner product as 
\begin{equation}
\label{innerproduct4example}
x\cdot y=\int_{[0,\infty)} e^{-\rho t}x(t)y(t)dt,
\end{equation}
and define the norm as $||x||=\sqrt{x\cdot x}$.  
Under this norm, we can show that $\mathcal{X}$ is complete, thus Assumption \ref{assumptionf_i_2} is satisfied.
\begin{lemma}
\label{hilbertspace}
The state space $\mathcal{X}$ relative to the mean-field game  (\ref{examplelqg}) is a Hilbert space.  
\end{lemma}
The proof can be found in Appendix A.

Third, under the inner product (\ref{innerproduct4example}), the objective function of the problem (\ref{examplelqg}) can be transformed to the form of (\ref{individualoptimization}):
\begin{equation}
\label{transformedcont}
\mathbb{E} \left( V_i(x_i,u_i)+F(\bar{x})\cdot x_i+G(\bar{x}) \right),
\end{equation}
where $V_i(x_i,u_i)=||x_i||^2+\int_0^{\infty} e^{-\rho t}ru_i(t)^2dt$, $F(\bar{x})=2c \bar{x}-2b \mathbb{I}$, $G(\bar{x})=||c \bar{x}-b \mathbb{I}||^2$, and $\mathbb{I}$ denotes the unit vector in $\mathcal{X}$.
Since $F(\bar{x})=2c \bar{x}-2b$ and $\nabla G(\bar{x})=2(c \bar{x}-b \mathbb{I} )$ are globally Lipschitz continuous, Assumption \ref{assumption_cost} is satisfied.  
Therefore, (\ref{examplelqg}) is a special case of the mean-field game (\ref{individualoptimization}).

\begin{remark}
In this example, the admissible control space is the set of  control strategies adapted to the filtration $\sigma(x_i(0),\pi_i(s),s\leq t)$. This represents the closed-loop perfect state information structure \cite{basar1995dynamic}, i.e. $u_i(t)$ depends on $\{x_i(s), 0\leq s \leq t\}$. We comment that the information structure of the mean-field game (\ref{individualoptimization}) can be  more general than closed-loop perfect state information structure. For instance, if $\bar{\mathcal{U}}_i$ only contains the control strategies that depend on the initial state $x_i(0)$, it represents the open-loop information structure. If each $u_i$ in $\bar{\mathcal{U}}_i$ satisfies that $u_i(t)$ depends on $x_i(t)$,  it corresponds to the close-loop feedback information structure. Such generalization is possible as long as the admissible control space $\bar{\mathcal{U}}_i$ remains a vector space.
\end{remark}

\subsubsection{Network Routing Games}
The third example concerns a splittable routing game \cite{cominetti2006network}, \cite{roughgarden2015local} on a network  $\mathcal{G}=(\mathcal{V}, \mathcal{E})$, where $\mathcal{V}$ is the set of vertices and $\mathcal{E}$ is the set of edges.  Assume there are $N$ players. Each of them (say, the $i$th agent) tries to allocate traffic $r_i$ to the paths from an origin $s_i\in \mathcal{V}$ to a destination $t_i\in \mathcal{V}$. This is done by choosing the flow $f_p^{(i)}$ such that $\sum_{p\in \mathcal{P}_i} f_p^{(i)}=r_i$, where $\mathcal{P}_i$ is the set of paths from $s_i$ to $t_i$, and $f_p^{(i)}$ is the rate allocated to path $p$. Given $f_p^{(i)}$, the total flow on a path, denoted as $f_p$, is the sum of all allocations to this path, i.e., $f_p=\sum_{i=1}^N f_p^{(i)}$. The flow on each edge $e\in \mathcal{E}$ is the total flow of all paths occupying this edge, i.e., $f_e=\sum_{p\in \mathcal{P}_i:e\in p} f_p$. 

Each edge $e\in \mathcal{E}$ of the network has a congestion cost that depends on the flow $f_e$. Denote this cost as $c_e(f_e)$, and assume it to be increasing, convex and continuously differentiable. The cost of each player is as follow: 
\begin{equation}
\label{networkroutinggame}
J_i(f)=\sum_{e\in p,p\in \mathcal{P}_i} c_e(f_e)f_p^{(i)}.
\end{equation}

Interestingly, the routing game (\ref{networkroutinggame}) is a special case of the proposed mean-field game (\ref{utilityfunction}). To see this, reformulate the decision of each agent as the vector $u_i=\{f_e^{(i)}, e\in \mathcal{E}\}$, where $f_e^{(i)}=\sum_{p\in \mathcal{P}_i: e\in p} f_p^{(i)}$ is the rate agent $i$ allocated to edge $e$. Clearly, $f_e=\sum_{i=1}^N f_e^{(i)}$.  Let $F(f)=\{F_e(f), e\in \mathcal{E}\}$ be a vector-valued function, and each  coordinate of $F(f)$ is $F_e(f)=c_e(f_e)$. Using this notation, the routing game  (\ref{networkroutinggame}) can be transformed to the following form:
\begin{equation}
\label{transformroutinggame}
J_i(u)=F(N\bar{u})\cdot u_i,
\end{equation}
where $\bar{u}=\dfrac{1}{N} \sum_{i=1}^N u_i$.
This is  a special case of (\ref{utilityfunction}) with $V_i=0$ and $G=0$. Note that the network imposes some constraint on the decision $u_i$. These are linear constraints, and we neglect them in the presentation.

\subsection{Objective of this Paper}
The objective of this paper is to study the connection between the mean-field game (\ref{individualoptimization}) and the social welfare optimization problem.  In most of the existing literature, a natural candidate for the social welfare is simply the sum of individual utilities. Unfortunately, it is shown that under this social welfare, the mean-field equilibrium is in general not efficient \cite{balandat2013efficiency}, and research effort has largely focused on characterizing and bounding the gap between the equilibrium solution and the efficient solution \cite{huang2007large}, \cite{yin2014efficiency}.

In this paper, instead of characterizing the gap, we ask the question of whether the mean-field equilibrium can be efficient for a modified social welfare. In other words, we would like to construct a social welfare optimization problem with some modified social welfare such that its optimal solution coincides with the mean-field equilibrium. 

In the rest of this paper, we will answer this question in two steps. First, we focus on the mean-field game  (\ref{individualoptimization}) and characterize its equilibrium as the solution to a set of mean-field equations. Second, we construct a social welfare optimization problem, and show that the solution to the mean-field equations coincides with the solution to the social welfare  optimization problem.

\section{Characterizing The Mean-Field Equilibrium}
This section  characterizes the equilibrium of the mean-field game (\ref{individualoptimization}) with a set of equations.  

There are several solution concepts for the game equilibrium, such as Nash equilibrium, Bayesian Nash equilibrium, dominant strategy equilibrium, among others. In mean-field games, we usually relax the Nash equilibrium solution concept by assuming that each agent is indifferent to an arbitrarily small change $\epsilon$. This solution concept is referred to as the $\epsilon$-Nash equilibrium, formally defined as follows:
\begin{definition}
\label{definitionofmeanfieldequilibrium}
$(u_1^*,\ldots,u_N^*)$ is an $\epsilon$-Nash equilibrium of the game (\ref{individualoptimization}) if  the following inequality holds 
\begin{equation}
 J_i(u_i^*,u_{-i}^*)\leq  J_i(u_i,u_{-i}^*)+\epsilon
\end{equation}
for all $i=1,\ldots,N$, and all $u_i\in \bar{\mathcal{U}}_i$, where $u_{-i}=(u_1,\ldots,u_{i-1},u_{i+1},\ldots,u_N)$ and $J_i(u_i,u_{-i})$ is the compact notation for (\ref{utilityfunction}) after plugging (\ref{systemdynamics}) in (\ref{utilityfunction}).
\end{definition}

At an $\epsilon$-Nash equilibrium, each agent can lower his cost by at most $\epsilon$ via deviating from the equilibrium strategy, given that all other players follow the equilibrium strategy. Therefore, the agents are motivated to play the equilibrium strategy if they are indifferent to a change of $\epsilon$ in their cost.

To characterize the $\epsilon$-Nash equilibrium, we note that the cost function (\ref{utilityfunction}) of the individual agent is only coupled through the mean-field terms $F(\bar{x})$ and $G(\bar{x})$. In the large population case, the impact of the control input for a single agent on the coupling term vanishes as the population size becomes large. Therefore, we can approximately treat the mean-field terms $F(\bar{x})$ and $G(\bar{x})$ as given, and each agent then faces an optimal response problem defined as follows:
\begin{align}
\label{optimalresponse}
\mu_i(y)\in \argmin_{u_i} &\mathbb{E} \left(V_i(x_i,u_i)+y\cdot  x_i \right) \\
\text{s.t. } 
&x_i=f_i(u_i, \pi_i) , \quad u_i\in \bar{\mathcal{U}}_i,
\end{align}  
where we use a deterministic value $y\in \mathcal{X}$ to replace $F(\bar{x})$. In (\ref{optimalresponse}), $\mu_i(y)$ denotes the optimal solution to the optimal response problem parameterized by $y$, and $G(\bar{x})$ is regarded as a constant in (\ref{optimalresponse}) that can be ignored. To avoid triviality, we impose the following assumption on (\ref{optimalresponse}):
\begin{assumption}
\label{assumptiononvandu}
For any $y\in \mathcal{X}$, the optimal response problem (\ref{optimalresponse}) admits at least one solution. 
\end{assumption}
Assumption \ref{assumptiononvandu} imposes some mild regularity conditions on the functional $V_i(x_i,u_i)$ and the admissible control set $\bar{\mathcal{U}}_i$. For instance, the solution to (\ref{optimalresponse}) exists if $V_i$ is continuous and $\bar{\mathcal{U}}_i$ is compact. However, the solution to (\ref{optimalresponse}) may also exist beyond these cases. Therefore, to maintain generality, we do not present the detailed technical condition for Assumption~\ref{assumptiononvandu} to hold. When the optimal response problem (\ref{optimalresponse}) have multiple solutions, $\mu_i(y)$ can be any one of these solutions.  

Based on this approximation,  $y$  generates a collection of agent responses. Ideally, the value $y$ should guide the individual agents to choose a collection of optimal responses $\mu_i(y)$ which, in return, collectively generate the mean-field term $F(\cdot)$ that is close to the approximation $y$. This suggests that we use the following equation systems to characterize the equilibrium of the mean-field game:
\begin{numcases}{\hspace{-0.7cm}}
   \mu_i(y)\in \argmin_{u_i\in \bar{\mathcal{U}}_i} \mathbb{E} \left(V_i(f_i(u_i, \pi_i),u_i)+y\cdot  f_i(u_i, \pi_i) \right)      \label{equationsystem1} 
   \\
   x_i=f_i(\mu_i(y), \pi_i)  \label{equationsystem2} \\
   y=F\left(\dfrac{1}{N}\sum_{i=1}^N \mathbb{E} x_i\right),  \label{equationsystem3}
\end{numcases}
where the mean-field term $y$ induces a collection of responses that generate $y$. In the rest of this subsection we show that the solution to this equation system is an $\epsilon$-Nash equilibrium of the mean-field game (\ref{individualoptimization}), and $\epsilon$ can be arbitrarily close to $0$ for sufficiently large $N$. For this purpose, we first prove the following two lemmas. These two lemmas are mainly used to set up the stage for the main result that will be introduced later. 
\begin{lemma}
\label{lemma1}
If there exists $C>0$ such that $\mathbb{E}||x_i||^2\leq C$ for all $i=1,\ldots,N$, and $F(\cdot)$ is globally Lipschitz continuous, then the following relation holds for each agent $i$:
\begin{equation}
\label{prooflemma}
\big|\mathbb{E} \big( F(\bar{x})\cdot x_i\big)-F(\mathbb{E}\bar{x})\cdot \mathbb{E} x_i \big|\leq \epsilon,
\end{equation}
where $\bar{x}\triangleq\dfrac{1}{N}\sum_{i=1}^N x_i$ and $0<\epsilon=O(\dfrac{1}{\sqrt{N}})$. 
\end{lemma}
See Appendix B for proof

The other lemma shows that removing the decision of a single agent does not significantly affect the value of $\mathbb{E}G(\cdot)$:
\begin{lemma}
\label{lemma2}
Define $\bar{x}\triangleq\dfrac{1}{N} \sum_{i=1}^Nx_i$, and let $\bar{x}_{-i}\triangleq\dfrac{1}{N} \sum_{j\neq i} x_j$, where $x_i$ denotes the state trajectory corresponding to $u_i$ in $\bar{\mathcal{U}_i}$, then we have the following relation:
\begin{equation}
\big |\mathbb{E}G(\bar{x})-\mathbb{E}G(\bar{x}_{-i})\big|\leq \epsilon,
\end{equation}
and $0<\epsilon=O(\dfrac{1}{N})$.
\end{lemma}
See Appendix C for proof.

Using the results of Lemma \ref{lemma1} and Lemma \ref{lemma2}, we can show that the solution to  the equation system (\ref{equationsystem1})-(\ref{equationsystem3}) is an $\epsilon$-Nash equilibrium of the mean-field game (\ref{individualoptimization}). 
\begin{theorem}
\label{theorem1}
The solution to the equation system (\ref{equationsystem1})-(\ref{equationsystem3}), if exists, is an $\epsilon_N$-Nash equilibrium of the mean-field game (\ref{individualoptimization}), and  $0<\epsilon=O(\dfrac{1}{\sqrt{N}})$. 
\end{theorem}
\vspace{-0.5cm}
\begin{pf*}{Proof.}
For notation convenience, we denote the solution to the equation system (\ref{equationsystem1})-(\ref{equationsystem3}) as $u_i^*$, $x_i^*$ and $y^*$, where $x_i^*$ is the state trajectory corresponding to $u_i^*$. 
According to Definition \ref{definitionofmeanfieldequilibrium}, to prove this theorem, we need to show that:
\begin{align}
\label{condition4meanfieldequilibrium}
&\mathbb{E} \bigg( V_i(x_i^*,u_i^*)+F(\bar{x}^*)\cdot x_i^*+G(\bar{x}^*)\bigg) \leq \epsilon+\nonumber\\ \mathbb{E} \bigg(&V_i(x_i,u_i)+F\bigg(\dfrac{1}{N}x_i+\bar{x}_{-i}^* \bigg)\cdot x_i+G\bigg(\dfrac{1}{N}x_i+\bar{x}_{-i}^*\bigg)\bigg)
\end{align}
for all $u_i\in \bar{\mathcal{U}}_i$, where $\epsilon=O(\dfrac{1}{\sqrt{N}})$, $\bar{x}^*\triangleq\dfrac{1}{N}\sum_{i=1}^N x_i^*$, $\bar{x}_{-i}^*\triangleq\dfrac{1}{N}\sum_{j\neq i} x_j^*$, and $x_i$ is the state trajectory corresponding to $u_i$. 
Based on Lemma \ref{lemma1}, we have 
\begin{align}
\label{resultoflemma1}
F(\mathbb{E}\bar{x})\cdot \mathbb{E} x_i - O(\dfrac{1}{\sqrt{N}}) \leq  \mathbb{E} \big( F(\bar{x})\cdot x_i\big) \leq \nonumber \\ F(\mathbb{E}\bar{x})\cdot \mathbb{E} x_i + O(\dfrac{1}{\sqrt{N}}).
\end{align}
Based on Lemma \ref{lemma2}, we have 
\begin{equation}
\label{resultoflemma2}
\mathbb{E}G(\bar{x}_{-i})- O(\dfrac{1}{\sqrt{N}}) \leq \mathbb{E}G(\bar{x}) \leq \mathbb{E}G(\bar{x}_{-i})+ O(\dfrac{1}{\sqrt{N}}).
\end{equation}
Apply (\ref{resultoflemma1}) and (\ref{resultoflemma2}) to the left-hard side of (\ref{condition4meanfieldequilibrium}), then the left-hand side of (\ref{condition4meanfieldequilibrium}) is upper bounded as follows:
\begin{align}
\label{intermediateinequality1}
\mathbb{E} \bigg( V_i(x_i^*,&u_i^*)+F(\bar{x}^*)\cdot x_i^*+G(\bar{x}^*)\bigg) \leq O(\dfrac{1}{\sqrt{N}})+\nonumber\\
 &\mathbb{E}V_i(x_i^*,u_i^*)+ F(\mathbb{E}\bar{x}^*)\cdot \mathbb{E} x_i^* +\mathbb{E} G(\bar{x}_{-i}^*).
\end{align}
Applying (\ref{resultoflemma1}) and (\ref{resultoflemma2}) to the right-hard side of (\ref{condition4meanfieldequilibrium}),  the right-hand side of (\ref{condition4meanfieldequilibrium}) is then lower bounded as follows:
\begin{align}
\label{intermediateinequality2}
\mathbb{E} \bigg(&V_i(x_i,u_i)+F\bigg(\dfrac{1}{N}x_i+\bar{x}_{-i}^* \bigg)\cdot x_i+G\bigg(\dfrac{1}{N}x_i+\bar{x}_{-i}^*\bigg)\bigg)  \nonumber\\ 
\geq &\mathbb{E}  V_i(x_i,u_i)+F\bigg(\dfrac{1}{N}\mathbb{E}\big(x_i+\sum_{j\neq i}x_j^* \big)\bigg)\cdot \mathbb{E}x_i +\nonumber \\
& \mathbb{E} G(\bar{x}_{-i}^*)-O(\dfrac{1}{\sqrt{N}}).
\end{align}

Based on (\ref{intermediateinequality1}) and (\ref{intermediateinequality2}), to prove that (\ref{condition4meanfieldequilibrium}) holds, it suffices to show that
\begin{align}
\label{provecondition}
\mathbb{E} V_i(x_i^*,&u_i^*)+F\bigg(\dfrac{1}{N}\mathbb{E}\sum_{i=1}^N x_i^*\bigg)\cdot \mathbb{E} x_i^* \leq O(\dfrac{1}{\sqrt{N}})+ \nonumber\\ 
&\mathbb{E}V_i(x_i,u_i)+F\bigg(\dfrac{1}{N}\mathbb{E}\big(x_i+\sum_{j\neq i}x_j^* \big)\bigg)\cdot \mathbb{E}x_i.
\end{align}
Since $||\mathbb{E} x_i||$ is bounded (see proof for Lemma \ref{lemma1}) and $F(\cdot)$ is Lipschitz continuous with  constant $L\geq 0$, we have:
\begin{align}
\label{provecondition_temp}
\bigg|F\bigg(\dfrac{1}{N}\mathbb{E}(x_i+&\sum_{j\neq i}x_j^* )\bigg)\cdot \mathbb{E}x_i-F\bigg(\dfrac{1}{N}\mathbb{E}\sum_{i=1}^N x_i^*\bigg)\cdot \mathbb{E}x_i\bigg|\leq \nonumber \\
&\big\Vert \dfrac{1}{N}L(\mathbb{E}x_i-\mathbb{E}x_i^*)\big\Vert \big\Vert \mathbb{E}x_i \big\Vert=O(\dfrac{1}{N}).
\end{align}
Therefore, combining (\ref{provecondition}) and (\ref{provecondition_temp}), it suffices to show that:
\begin{align*}
\mathbb{E} V_i(x_i^*,&u_i^*)+F\bigg(\dfrac{1}{N}\mathbb{E}\sum_{i=1}^N x_i^*\bigg)\cdot \mathbb{E} x_i^* \leq \nonumber\\ &\mathbb{E}V_i(x_i,u_i)+F\bigg(\dfrac{1}{N}\mathbb{E}\sum_{i=1}^N x_i^*\bigg)\cdot \mathbb{E}x_i+O(\dfrac{1}{\sqrt{N}}).
\end{align*}
Note that based on (\ref{equationsystem3}), $F\bigg(\dfrac{1}{N}\mathbb{E}\sum_{i=1}^N x_i^*\bigg)=y^*$, which is equivalent to:
\begin{equation*}
\mathbb{E} V_i(x_i^*,u_i^*)+y^*\cdot \mathbb{E} x_i^* \leq \mathbb{E}V_i(x_i,u_i)+y^*\cdot \mathbb{E}x_i+O(\dfrac{1}{\sqrt{N}}).
\end{equation*}
This obviously holds based on (\ref{equationsystem1}), which completes the proof.
\end{pf*}
Theorem \ref{theorem1} indicates that each agent is motivated to follow the equilibrium strategy $u_i^*$ as deviating from this strategy can only decrease the individual cost by a negligible amount $\epsilon$. Furthermore, this $\epsilon$ can be arbitrarily small, if the population size is sufficiently large. Note that the mean-field equation system (\ref{equationsystem1})-(\ref{equationsystem3}) is not the unique way to characterize the $\epsilon$-Nash equilibrium of the mean-field game (\ref{individualoptimization}). The game may have other $\epsilon$-Nash equilibria with different values of $\epsilon$. However, in this paper, we only focus on the mean-field equations (\ref{equationsystem1})-(\ref{equationsystem3}). In the rest of this paper, the mean-field equilibrium of the game (\ref{individualoptimization}) always refers to the solution to the mean-field equations (\ref{equationsystem1})-(\ref{equationsystem3}).


\begin{remark} 
\label{comparison}
It is interesting to discuss the connections and differences between the mean field equations in our paper and those in \cite{huang2007large}. 
The deference between these two results is clear: we consider a set of mean-field equations for a finite population of agents, while in \cite{huang2007large} the authors used a continuum model to approximate the finite game. This leads to slightly different forms of the mean-field equations: we use empirical mean $\dfrac{1}{N}\sum_{i=1}^N \mathbb{E} x_i$ to represent the mean-field, while  \cite{huang2007large} assumes a distribution over the agent parameter and uses integral of system state over this distribution to represent the mean (see (4.8) in \cite{huang2007large}). On the other hand,  these two models are connected if we let $N$ tends to infinity: if the empirical mean converges to the integral, then the mean-field equations  (\ref{equationsystem1})-(\ref{equationsystem3}) are equivalent to those in \cite{huang2007large}. However, we emphasize that this is only true if the empirical mean does converge as $N$ goes to infinity. Otherwise, the mean-field equations (\ref{equationsystem1})-(\ref{equationsystem3}) are not well-defined as $N$ approaches infinity. Therefore, we can regard the mean-field equations in \cite{huang2007large} as the limiting case of our equations (\ref{equationsystem1})-(\ref{equationsystem3}) under additional assumptions that guarantee the convergence of the empirical mean $\dfrac{1}{N}\sum_{i=1}^N \mathbb{E} x_i$. Such assumption was imposed in \cite{huang2007large} (see Assumption H3 in \cite{huang2007large}). We do not have these assumptions in our paper,  since we directly derive mean-field equations for a finite population of agents. 
\end{remark}


\section{Connection to Social Welfare Optimization}
This section focuses on the connection between the mean-field game (\ref{individualoptimization}) and the social welfare optimization problem. 
Such connection is typically referred to as ``efficiency": we say that the mean-field equilibrium is efficient if it maximizes the social welfare.  In the literature, some attempts have been made to draw connections between the mean-field game and the social welfare optimization problem. Most of these works consider the social welfare optimization problem to be maximizing the total utility (or equivalently, minimizing the total cost) of all agents, which in our context can be formulated as follows:
\begin{align}
\label{socialwelfareoptbenchmark}
\min_{(u_1,\ldots,u_N)} &\sum_{i=1}^N \mathbb{E} \left(V_i(x_i,u_i)+F(\bar{x})\cdot  x_i +G(\bar{x})\right) \\
\text{s.t. } 
\label{abstractdynamic}
&x_i=f_i(u_i, \pi_i), \quad u_i\in \bar{\mathcal{U}}_i, \quad \forall i\in \mathcal{N}. 
\end{align}
where $\mathcal{N}=\{1,\ldots,N\}$.
Since the cost function (\ref{socialwelfareoptbenchmark}) represents the total cost of the entire population, from the efficiency point of view, it is desirable to have the mean-field equilibrium to be the optimal solution to (\ref{socialwelfareoptbenchmark}). However, it is shown that this statement is not true in general \cite{huang2007large}, \cite{yin2014efficiency}, \cite{balandat2013efficiency}. Therefore, many  existing works along this line focused on characterizing the gap between the mean-field equilibrium and the optimal solution to (\ref{socialwelfareoptbenchmark}).

Different from these works, we construct a modified  social welfare optimization problem so that the mean-field equilibrium achieves exact social optima. This can be done by introducing a virtual agent in the system with a cost function $\phi:\mathcal{X}\rightarrow \mathbb{R}$, and consider the following social welfare optimization problem that includes this virtual cost:
\begin{align}
\label{socialwelfare_constructed}
&\min_{u_1,\ldots,u_N,z} \mathbb{E} \left(\sum_{i=1}^N V_i(x_i,u_i)+\phi(z)\right)\\
\text{s.t.} &
\begin{cases}
z= \dfrac{1}{N}   \sum_{i=1}^N \mathbb{E} x_i, \\
x_i=f_i(u_i, \pi_i), \quad u_i\in \bar{\mathcal{U}}_i, \quad \forall i\in \mathcal{N}. 
\end{cases}
\end{align}
where $z$ is the decision of the virtual agent. Compared to the classical social welfare optimization problem (\ref{socialwelfareoptbenchmark}), the constructed problem has an augmented decision variable, and introduces an additional constraint $z= \dfrac{1}{N}   \sum_{i=1}^N \mathbb{E} x_i$. This constraint is inspired by the supply-demand model in microeconomics, where the virtual agent acts as a single supplier, and all other agents are the demands. This constraint requires that the supply and the demands are balanced.

In the remainder of the section, we will establish conditions under which we can draw connections between the mean-field equilibrium and the solution to the social welfare optimization  problem (\ref{socialwelfare_constructed}). Finding such conditions is useful from at least two perspectives. First, if we are given a mean-field game, we can construct the social welfare optimization problem  to evaluate the efficiency of the mean-field equilibrium. Second, if we are given a social welfare optimization problem (\ref{socialwelfareoptbenchmark}), then we can design the mean-field game  so as to control the population to operate at the socially optimal point.

\subsection{Connections under Strong Duality}
This subsection shows that the solution to the social welfare optimization problem (\ref{socialwelfare_constructed}) is a mean-field equilibrium to the game (\ref{individualoptimization}) if  (\ref{socialwelfare_constructed}) has strong duality.

For this purpose, we first introduce the concept of strong duality for the social welfare optimization problem (\ref{socialwelfare_constructed}). With slight abuse of notation, we drop the dependence of the objective function of (\ref{socialwelfare_constructed}) on $x$, and compactly denote  (\ref{socialwelfare_constructed}) as follows:
\begin{align}
\label{compactsocialwelfare}
P^*=&\min_{u,z} J_s(u,z)\\
\text{s.t.}&
\begin{cases}
g(u,z)=0 \\
z\in \mathcal{X}, u_i\in \bar{\mathcal{U}}_i,\quad \forall i=1,\ldots,N,
\end{cases}
\end{align}
where $P^*$ is the optimal value of the social welfare optimization problem (\ref{compactsocialwelfare}), $u=(u_1,\ldots,u_N)$ is the vector of control inputs, $J_s(u,z)=\mathbb{E} \left(\sum_{i=1}^N V_i(f_i(u_i,\pi_i),u_i)+\phi(z)\right)$,  and  $g(u,z)=  \mathbb{E}\sum_{i=1}^N f_i(u_i,\pi_i)-Nz$. Using this notation, the  Lagrangian of problem (\ref{compactsocialwelfare}) can be defined as follows:
\begin{align}
\label{lagrangian}
L(u,z,\lambda)=J_s(u,z)+\lambda\cdot g(u,z).
\end{align}
where $\lambda\in \mathcal{X}$ is the Lagrange multiplier for the constraint $g(u,z)=0$.

\begin{remark}
\label{completeness}
In a more general setting, the Lagrange multiplier is in the dual space of $\mathcal{X}$, and the inner product term in (\ref{lagrangian}) should be replaced with a bounded linear operator evaluated at point $g(u,z)$ \cite[Chap. 8]{luenberger1997optimization}. However, in our problem, since the state space $\mathcal{X}$ is a Hilbert space, we can select the dual space of $\mathcal{X}$ to be itself, and the bounded linear operator reduces to the inner product on $\mathcal{X}$. Therefore, the expression (\ref{lagrangian})  relies on the fact that $\mathcal{X}$ is a Hilbert space. When $\mathcal{X}$ is not a complete space, the mean-field equations can not be connected to the social welfare optimization problem. 
\end{remark}

Given the Lagrangian (\ref{lagrangian}), we first treat the multiplier as given and define the mapping $D:\mathcal{X}\rightarrow \mathbb{R}$:
\begin{equation}
\label{intermediadual}
D(\lambda)=\inf_{u_i\in \bar{\mathcal{U}}_i, z\in \mathcal{X}} L(u,z,\lambda),
\end{equation}
then the dual problem of the social welfare optimization problem (\ref{socialwelfare_constructed}) is defined as follows:
\begin{equation}  
\label{dualprolem}
D^*=\max_{\lambda\in \mathcal{X}} D(\lambda)
\end{equation}
where $D^*$ is the optimal value of the dual problem. 
When the dual problem (\ref{dualprolem}) admits a solution, and the optimal value of the dual problem (\ref{dualprolem}) coincides with that of the primal problem (\ref{compactsocialwelfare}), then we say the optimization problem (\ref{socialwelfare_constructed}) has strong duality. Formally, we define it as follows:
\begin{definition}
\label{strongdualitydefition}
The optimization problem (\ref{socialwelfare_constructed}) has strong duality if $P^*=D^*$ and there exists $\lambda^*\in \mathcal{X}$ such that $D^*=D(\lambda^*)$.
\end{definition}
Note that the definition of strong duality not only requires the duality gap to be zero, but also requires the dual problem to have a finite solution $\lambda^*$. This is slightly stronger than only requiring zero duality gap between the primal problem (\ref{compactsocialwelfare}) and the dual problem (\ref{dualprolem}). In general, the existence of a finite multiplier to (\ref{dualprolem}) can be easily guaranteed under mild constraint qualifications (e.g., Slater's condition) \cite[Chap. 8]{luenberger1997optimization}. In our paper, when we say that the social welfare optimization problem has strong duality, it indicates the problem already satisfies certain constraint qualifications so that the solution to (\ref{dualprolem}) exists.

Under strong duality, we can establish connections between the mean-field equilibrium and the social welfare optimization problem.  These connections are summarized in the next a few theorems and corollaries, which are the main results of this paper. 
\begin{theorem}
\label{mainresult1}
Let $\phi:\mathcal{X} \rightarrow \mathbb{R}$ be a Fr\'echet differentiable functional such that $\nabla \phi(z)=NF(z)$, $\forall z\in \mathcal{X}$. Assume that the social welfare optimization problem (\ref{socialwelfare_constructed}) has strong duality, then any socially optimal solution to (\ref{socialwelfare_constructed}) is a mean-field equilibrium to the game (\ref{individualoptimization}). 
\end{theorem}
\vspace{-0.5cm}
\begin{pf*}{Proof.}
Since the social welfare optimization   (\ref{socialwelfare_constructed}) has strong duality, then there exists $\lambda^*$ such that $P^*=D^*=D(\lambda^*)$. Note that due to weak duality, this indicates that $\lambda^*$ is the optimal solution to the dual problem (\ref{dualprolem}), i.e., $D^*=\inf_{u_1\in \bar{\mathcal{U}}_1,\ldots, u_N\in \bar{\mathcal{U}}_N, z\in \mathcal{X}} L(u,z,\lambda^*)$. Let $(u^*,z^*)$ be the optimal solution to (\ref{socialwelfare_constructed}), then $(u^*,z^*)$ satisfies the constraint $z^*=\dfrac{1}{N}\sum_{i=1}^N \mathbb{E} f(u_i^*,\pi_i)$, and we have the following inequalities:
\begin{align}
\label{proofdetailtemp}
D^*&=\inf_{u_1\in \bar{\mathcal{U}}_1,\ldots, u_N\in \bar{\mathcal{U}}_N, z\in \mathcal{X}} L(u,z,\lambda^*) \nonumber \\
&\leq L(u^*,z^*,\lambda^*)=J_s(u^*,z^*)+\lambda^* \cdot g(u^*,z^*) \nonumber \\
&=J_s(u^*,z^*)=P^*.
\end{align} 
Due to strong duality, $P^*=D^*$. Therefore, equality holds in (\ref{proofdetailtemp}), indicating that $(u^*,z^*)$ satisfies the following:
\begin{equation}
\label{temp1}
(u^*,z^*)\in \argmin_{u_1\in \bar{\mathcal{U}}_1,\ldots, u_N\in \bar{\mathcal{U}}_N, z\in \mathcal{X}} L(u,z,\lambda^*).
\end{equation}
Since $L$ can be decomposed in terms of $u^i$ and $z$, the relation (\ref{temp1}) is equivalent to the following:
\begin{numcases}{\hspace{-0.5cm}}
 u_i^*\in \argmin_{u_i\in \bar{\mathcal{U}}_i} \mathbb{E} \left(V_i(f_i(u_i, \pi_i),u_i)+\lambda^*\cdot  f_i(u_i, \pi_i) \right)      \label{tempequation1}
   \\
   z^*\in \argmin_{z\in \mathcal{X}} \phi(z)-N\lambda^* \cdot z \label{tempequation2} 
\end{numcases}
The first-order optimality condition of (\ref{tempequation2}) yields $\nabla \phi(z^*)=N\lambda^*$. Since $\nabla  \phi(z)=NF(z)$, we have $F(z^*)=\lambda^*$. Therefore, the above equation sets can be reduced to the following:
\begin{numcases}{\hspace{-0.5cm}}
   u_i^*\in \argmin_{u_i\in \bar{\mathcal{U}}_i} \mathbb{E} \left(V_i(f_i(u_i, \pi_i),u_i)+\lambda^*\cdot  f_i(u_i, \pi_i) \right)      \label{equationsystem_proof1}
   \\
   \lambda^*=F \left(\dfrac{1}{N}\sum_{i=1}^N \mathbb{E} f_i(u_i^*,\pi_i)\right)  \label{equationsystem_proof2} 
\end{numcases}
It can be verified that (\ref{equationsystem_proof1})-(\ref{equationsystem_proof2}) is equivalent to the mean-field equations (\ref{equationsystem1})-(\ref{equationsystem3}). Therefore, $(u_1^*,\ldots,u_N^*)$ is a mean-field equilibrium. This completes the proof.
\end{pf*}
Theorem \ref{mainresult1} shows  that any socially optimal solution is a mean-field equilibrium if $\nabla \phi(z)=NF(z)$. This indicates that  we can construct $\phi(\cdot)$ from $F(\cdot)$ to study the mean-field equilibrium. Checking the existence of $\phi(\cdot)$ and computing $\phi(\cdot)$ from $F(\cdot)$ are non-trivial for an arbitrary $F(\cdot)$ in general vector spaces. However, since the majority of the literature formulates  the cost function as an integral over time \cite{huang2007large}, \cite{yin2012synchronization}, $F(\cdot)$  typically has special structures: it represents a trajectory over time, and its value at each time $t$ only depends on $x(t)$.  In this case, we can easily derive $\phi(\cdot)$ using variational analysis. An example is given in Section VI-B.

Theorem  \ref{mainresult1}  only shows that any socially optimal solution is the mean-field equilibrium.  This does not necessarily mean that any mean-field equilibrium is also socially optimal.  However, the other direction of the relation also holds when the mean-field equations have a unique solution. This can be summarized in the following corollary:
\begin{corollary}
\label{corollary2}
Let $\phi:\mathcal{X}\rightarrow \mathbb{R}$ be a Fr\'echet differentiable functional such that $\nabla \phi(z)=NF(z)$, $\forall z\in \mathcal{X}$. Assume that the social welfare optimization problem (\ref{socialwelfare_constructed}) has strong duality, and the mean-field game (\ref{socialwelfare_constructed}) has a unique mean-field equilibrium, then the mean-field equilibrium to (\ref{individualoptimization})  is the globally optimal solution to the social welfare optimization problem (\ref{socialwelfare_constructed}). 
\end{corollary}
The proof of the corollary follows easily  from Theorem \ref{mainresult1}, and is therefore omitted. Corollary \ref{corollary2} can be used to check the efficiency of the mean-field equilibrium  when the mean-field equations admit at most one solution.

\subsection{Special Case with Monotone Mean-Field Coupling}
In general, the mean-field equations may admit multiple solutions. According to Theorem \ref{mainresult1}, the best mean-field equilibrium among these solutions is the optimal solution to the social welfare optimization problem, but there may exist other mean-field equilibria that are not socially optimal. In this subsection, we show  that this complication can be resolved if the following monotonicity condition is imposed on the mean-field coupling term $F(\cdot)$:
\begin{definition}[monotone mean-field coupling]
\label{monotonecondition}
The mean-field coupling term $F(x)$ is monotone with respect to $x\in \mathcal{X}$, if $(F(x)-F(x'))\cdot (x-x')\geq 0$ for any $x,x'\in \mathcal{X}$.
\end{definition}

Under this condition, $\phi(\cdot)$ is convex, and we can derive a stronger result than Theorem \ref{mainresult1}, where the relation between the mean-field equilibrium and the socially optimal solutions can go either way: 
\begin{theorem}
\label{mainresult1_2}
Let $\phi:\mathcal{X}\rightarrow \mathbb{R}$ be a Fr\'echet differentiable functional such that $\nabla \phi(z)=NF(z)$, $\forall z\in \mathcal{X}$. Assume that the social welfare optimization problem (\ref{socialwelfare_constructed}) has strong duality, and assume that $F(\cdot)$ is monotone, then $(u_1^*,\ldots,u_N^*)$ is the mean-field equilibrium to (\ref{socialwelfare_constructed}) if and only if it is the globally optimal solution to the social welfare optimization problem (\ref{socialwelfare_constructed}).
\end{theorem}
\vspace{-0.5cm}
\begin{pf*}{Proof.}
Based on Theorem \ref{mainresult1}, the socially optimal solution is a mean-field equilibrium.  Therefore, it suffices to show the other direction also holds. For notational convenience,  let $\bar{u}=(\bar{u}_1,\ldots,\bar{u}_N)$ be the solution to the mean-field equations (\ref{equationsystem1})-(\ref{equationsystem3}). Define $\bar{z}=\dfrac{1}{N}\mathbb{E}\sum_{i=1}^N f_i(\bar{u}_i,\pi_i)$, and let $\bar{y}=F(\bar{z})$. Since $F(z)=\dfrac{1}{N}\nabla \phi(z)$, we have $\nabla \phi(\bar{z})=N\bar{y}$. Since $F(z)$ is monotone, $\phi(z)$ is convex. Therefore, $\nabla \phi(\bar{z})=N\bar{y}$ indicates that:
\begin{equation}
\label{reverseproblem}
\bar{z}\in \argmin_{z\in \mathcal{X}}\phi(z)-N\bar{y}\cdot z.
\end{equation}
Due to (\ref{equationsystem1}), we also have:
\begin{equation}
\label{reverseproblem2}
\bar{u}_i\in \argmin_{u_i\in \bar{\mathcal{U}}_i} \mathbb{E} \left(V_i(f_i(u_i, \pi_i),u_i)+\bar{y}\cdot f_i(u_i, \pi_i) \right).
\end{equation}
The above two equations together indicate that $(\bar{u},\bar{z})$ is the optimal solution to the following optimization problem:
\begin{align}
\label{proofcondition}
&\min_{u,z} \sum_{i=1}^N \mathbb{E} V_i(x_i,u_i)+\phi(z)+\bar{y}\cdot (\sum_{i=1}^N \mathbb{E} x_i-Nz) \\
&\text{s.t.} 
\begin{cases}
\label{abstractdynamic}
x_i=f_i(u_i, \pi_i) \\
 u_i\in \bar{\mathcal{U}}_i, z\in \mathcal{X}.
\end{cases}
\end{align}
In other words, $(\bar{u},\bar{z})$ satisfies:
 \begin{equation}
\label{temp1234567}
(\bar{u},\bar{z})\in \argmin_{u_1\in \bar{\mathcal{U}}_1,\ldots, u_N\in \bar{\mathcal{U}}_N, z\in \mathcal{X}} L(u,z,\bar{y}).
\end{equation}
Note that due to weak duality, we have:
\begin{equation}
\label{oneside}
L(\bar{u},\bar{z},\bar{y})\leq D^*\leq P^*. 
\end{equation}
On the other hand, we also have:
\begin{align}
\label{theotherside}
L(\bar{u},\bar{z},\bar{y})&=J_s(\bar{u},\bar{z})+\bar{y}\cdot g(\bar{u},\bar{z}) \nonumber \\
&=J_s(\bar{u},\bar{z})\geq P^*,
\end{align}
where the last inequality is due to the fact that $P^*$ is the minimum value of $J_s(u,z)$ among all $(u,z)$ such that $z=\dfrac{1}{N} \mathbb{E}\sum_{i=1}^N f_i(u_i,\pi_i)$, and $(\bar{u},\bar{z})$ is one of them. Combing (\ref{oneside}) and (\ref{theotherside}), we have $L(\bar{u},\bar{z},\bar{y})=P^*$, thus $(\bar{u},\bar{z})$ is the globally optimal solution to the social welfare optimization problem (\ref{socialwelfare_constructed}). This completes the proof.
\end{pf*}

This theorem establishes equivalence between the mean-field equilibrium and the socially optimal solutions. When the mean-field term is monotone, the solution set of the mean-field equations is the same as that of the social welfare optimization problem as long as (\ref{socialwelfare_constructed}) has strong duality. Regarding this result, an  interesting special case is where the social welfare optimization problem is convex. To ensure convexity, we introduce the following conditions:
\begin{assumption}
\label{convexityassumption}
(i) $f_i(u_i,\pi_i)$ is affine with respect to $u_i$, $\forall \pi_i\in \Pi$,  
(ii) $V_i(x_i,u_i)$ is convex with respect to $(x_i,u_i)$, 
(iii) $\mathcal{X}_i$ and $\mathcal{U}_i$ are convex, 
(iv) $F(\cdot)$ is monotone.
\end{assumption}
It can be easily verified that the social welfare optimization problem is convex with respect to $(x_1,\ldots,x_N,z)$ under  Assumption \ref{convexityassumption}. In this case, the social welfare optimization problem has strong duality  under mild constraint qualifications. Therefore, we have the following corollary: 
\begin{corollary}
\label{corollary1}
Let $\phi:\mathcal{X}\rightarrow \mathbb{R}$ be a Fr\'echet differentiable functional such that $\nabla \phi(z)=NF(z)$, $\forall z\in \mathcal{X}$. Under Assumption \ref{convexityassumption}, $(u_1^*,\ldots,u_N^*)$ is the mean-field equilibrium to (\ref{individualoptimization}) if and only if it is the globally optimal solution to the social welfare optimization problem (\ref{socialwelfare_constructed}). 
\end{corollary}
This corollary can be easily proved based on Theorem~\ref{mainresult1_2}:  Assumption \ref{convexityassumption} guarantees that the social welfare optimization problem (\ref{socialwelfare_constructed}) is convex with respect to  $(x_1,\ldots,x_N,z)$. Therefore, based on duality theory \cite[p. 224]{luenberger1997optimization}, the social welfare optimization problem (\ref{socialwelfare_constructed}) has strong duality. Furthermore,  according to Theorem  \ref{mainresult1_2},  the solution set of the mean-field equations is the same as that of  the social welfare optimization problem.

In the more general case, Assumption \ref{convexityassumption} may not be satisfied, and the social welfare optimization problem may be non-convex with  respect to $(x_1,\ldots,x_N,z)$. Therefore, to apply the result of the theorems, we need to check whether the social welfare optimization problem has strong duality, especially when Assumption \ref{convexityassumption} is not satisfied.  Along this direction, many sufficient conditions have been developed to check the strong duality for non-convex optimization \cite{flores2013strong}, \cite{carcamo2016strong}, \cite{flores2012complete}. 
Due to space limit, we will not present the technical details of these works. Instead, we will provide an  example in Section VI, where the social welfare optimization problem does not satisfy Assumption \ref{convexityassumption}, but Theorem \ref{mainresult1} and Theorem \ref{mainresult1_2} can  still be applied.

\begin{remark}
\label{remarkoncontrolmean}
In this paper, we formulate the mean-field term $F(\cdot)$ to depend on the average of the population {\bf state}. However, the proposed method still works when the mean-field is the average of the {\bf control} decisions. In this case, the individual cost functional (\ref{utilityfunction}) is defined as $J_i(x_i,u_i,F(m), G(m))=V_i(x_i,u_i)+F\bigg(\dfrac{1}{N}\sum_{i=1}^N u_i\bigg)\cdot u_i+G\bigg(\dfrac{1}{N}\sum_{i=1}^N u_i\bigg)$. Similar result can be obtained using the same approach.
\end{remark}

\subsection{Relation to Potential Games}
The proposed  social welfare optimization is different from the potential function in potential games. This subsection discusses their differences.  

First, we show that in general, the mean-field game (\ref{individualoptimization}) is not a potential game, thus  (\ref{socialwelfare_constructed}) is not a potential function.  
To construct a counter-example, consider a group of agents with the following cost function:
\begin{equation}
\label{examplenopotential}
J_i(x_1,\cdots,x_N)= (x_i-1)^2+ x_i \mbox{log}\bar{x} ,
\end{equation}
Note that on the one hand, it is easy to verify that (\ref{examplenopotential}) satisfies Assumption \ref{convexityassumption}. Therefore, based on Theorem \ref{mainresult1_2},  the mean-field equilibrium is equivalent to the modified social welfare optimization problem (\ref{socialwelfare_constructed}). On the other hand, based on Theorem 4.5 in \cite{monderer1996potential},  (\ref{examplenopotential})  is a potential game if and only if: 
\begin{equation}
\label{ifonlyifpotential}
\dfrac{\partial^2 J_i}{\partial x_i \partial x_j} = \dfrac{\partial^2 J_j}{\partial x_i \partial x_j}.
\end{equation}
It can be verified that (\ref{ifonlyifpotential}) does not hold for (\ref{examplenopotential}). Therefore, it is not a potential game.

Second, under additional assumptions (i.e., $F(\bar{x})=c\sum_{i=1}^N x_i$), the mean-field game may reduce to a potential game. In this case, we can either solve the potential function minimization or the social welfare optimization problem to derive the mean-field equilibrium. However, the solutions to these two methods are different: the potential function method provides a Nash equilibrium of the game, while the social welfare optimization leads to the $\epsilon$-Nash equilibrium. This is a relaxed solution concept.
On the other hand, from the computational perspective, it is more attractive to use the social welfare optimization instead of the potential function. This is because the computational complexity of solving the potential minimization problem increases as the number of agents increases, while the social welfare optimization enables a decentralized scheme where the computation time is irrelevant with respect to the number of agents in the game \cite{ma2013decentralized}. This property is important in large-scale game problems.  More details on the computation of mean-field equilibria can be found in the next section.

\section{Implication on Existence, Uniqueness and Computation}
Our results have some interesting implications on the existence, uniqueness and computation of the mean-field equilibrium. These implications extend the results in the literature to more general cases. We discuss these implications in this section.

\subsection{Existence of the Mean-Field Equilibrium}
The existence and uniqueness of the mean-field equilibrium is a problem of fundamental importance in mean-field games. This problem has been extensively studied in the literature \cite{lasry2007mean},  \cite{ahuja2016wellposedness}, \cite{carmona2013probabilistic}, \cite{carmona2015probabilistic}, \cite{cardaliaguet2010notes}, and many of these works are based on fixed point analysis. In this section, we provide a novel approach inspired by the connection between the mean-field game and the social welfare optimization problem. Our approach extends existing results to more general cases.   We start with the following result:
\begin{proposition}
\label{existencetheorem}
There exists a mean-field equilibrium to (\ref{individualoptimization}) if the social welfare optimization problem (\ref{socialwelfare_constructed}) has strong duality, where  $\phi(\cdot)$ is such that $\nabla \phi(z)=NF(z)$, $\forall z\in \mathcal{X}$. 
\end{proposition}
The proof  directly follows from Theorem~\ref{mainresult1}: since the optimal solution to (\ref{socialwelfare_constructed}) is a mean-field equilibrium, if (\ref{socialwelfare_constructed}) admits a solution, then the mean-field equations also admit a solution. An interesting fact about Theorem~\ref{existencetheorem} is that it draws connections between the existence of the mean-field equilibrium and the strong duality of (\ref{socialwelfare_constructed}). This enables us to check the existence of the mean-field equilibrium by verifying the strong duality of (\ref{socialwelfare_constructed}). 


Proposition \ref{existencetheorem}  provides a sufficient condition for the existence of the mean-field equilibrium. In fact, we can show that this condition is also necessary under additional assumptions on the mean-field coupling term $F(\cdot)$. This is summarized in the following theorem: 
\begin{proposition}
\label{existenceNAScondition}
Assume that $F(\cdot)$ is monotone, then there exists a mean-field equilibrium to the game (\ref{individualoptimization})  if and only if the social welfare optimization problem (\ref{socialwelfare_constructed}) has strong duality, where $\phi(\cdot)$ is such that $\nabla \phi(z)=NF(z)$, $\forall z\in \mathcal{X}$. 
\end{proposition}
The proof  is a byproduct of the proof of Theorem \ref{mainresult1_2}, which directly follows from (\ref{oneside}) and (\ref{theotherside}).  We comment that most related results in existing works only have sufficient conditions for the existence of the mean-field equilibrium \cite{lasry2007mean}, \cite{ahuja2016wellposedness}, \cite{carmona2013probabilistic}, \cite{carmona2015probabilistic}, \cite{cardaliaguet2010notes}. In view of this, the significance of  Theorem \ref{existenceNAScondition} is that it provides an existence condition that is both necessary and sufficient.  Using this result, we can not only show that the mean-field equilibrium exists under strong duality,  but also show that the mean-field equilibrium does not exist when strong duality does not hold.

\subsection{Uniqueness of the Mean-Field Equilibrium}
We has the following result on the uniqueness of the mean-field equilibrium:
\begin{proposition}
\label{uniquenesstheorem}
There is a unique mean-field equilibrium to the game (\ref{individualoptimization}) if the social welfare optimization problem (\ref{socialwelfare_constructed}) is strictly convex with respect to $(u_1,\ldots,u_N,z)$, where $\phi(\cdot)$ is such that $\nabla \phi(z)=NF(z)$, $\forall z\in \mathcal{X}$.
\end{proposition}
The proof  follows  from Corollary~\ref{corollary1}: under the assumptions, the solution set to the mean-field equations is equivalent to that of the social welfare optimization problem. As (\ref{socialwelfare_constructed}) is strictly convex, it has a unique solution. Therefore, the mean-field equations also have a unique solution. 

Proposition  \ref{uniquenesstheorem} has interesting connections to many existing works on the uniqueness of the mean-field equilibrium. Although most of these works focus on continuum mass model, if we consider the finite counterpart of these works and adapt their models to our context, then we can roughly divide the conditions in these works in three categories. First, the cost function $V_i(f(u_i,\pi_i),u_i)$ is assumed to be strictly convex with respect to $u_i$, and the coupling term $F(\cdot)$ is monotone \cite{ahuja2016wellposedness},  \cite{carmona2013probabilistic}, \cite{carmona2015probabilistic}. Second, the cost function $V_i(f(u_i,\pi_i),u_i)$ is convex with respect to $u_i$, and the coupling term $F(\cdot)$ is strictly increasing \cite{gomes2013extended}, \cite{cardaliaguet2010notes}. Third, the cost function is at least convex, and the coupling term is at least monotone, but either one of them holds strictly \cite{lasry2007mean}. We note that these conditions can be all recovered by Theorem {\ref{uniquenesstheorem}:  it is not hard to verify that all these conditions essentially ensure the social welfare optimization problem (\ref{socialwelfare_constructed}) to be strictly convex with respect to $(u_1,\ldots,u_N,z)$. According to Theorem {\ref{uniquenesstheorem}}, there is a unique mean filed equilibrium for the mean-field game. Therefore, Theorem~{\ref{uniquenesstheorem}} connects the uniqueness of the mean-field equilibrium to the strict convexity of an optimization problem, providing a novel interpretation of these existing results. 

To summarize, we have  established the connections between the mean-field equilibrium and the social welfare optimization problem, and discussed the implications on the existence and uniqueness of the mean-field equilibrium. To better understand these results, we graphically summarize these connections using a diagram in Figure~\ref{diagram}. In this figure,  S.O. stands for the socially optima solution to (\ref{socialwelfare_constructed}), MFE stands for the mean-field equilibrium, the arrow indicates the conclusion we can draw, and the text on the arrow denotes the conditions for the conclusion.

\begin{figure}[bt]%
\centering
\includegraphics[width = 0.85\linewidth]{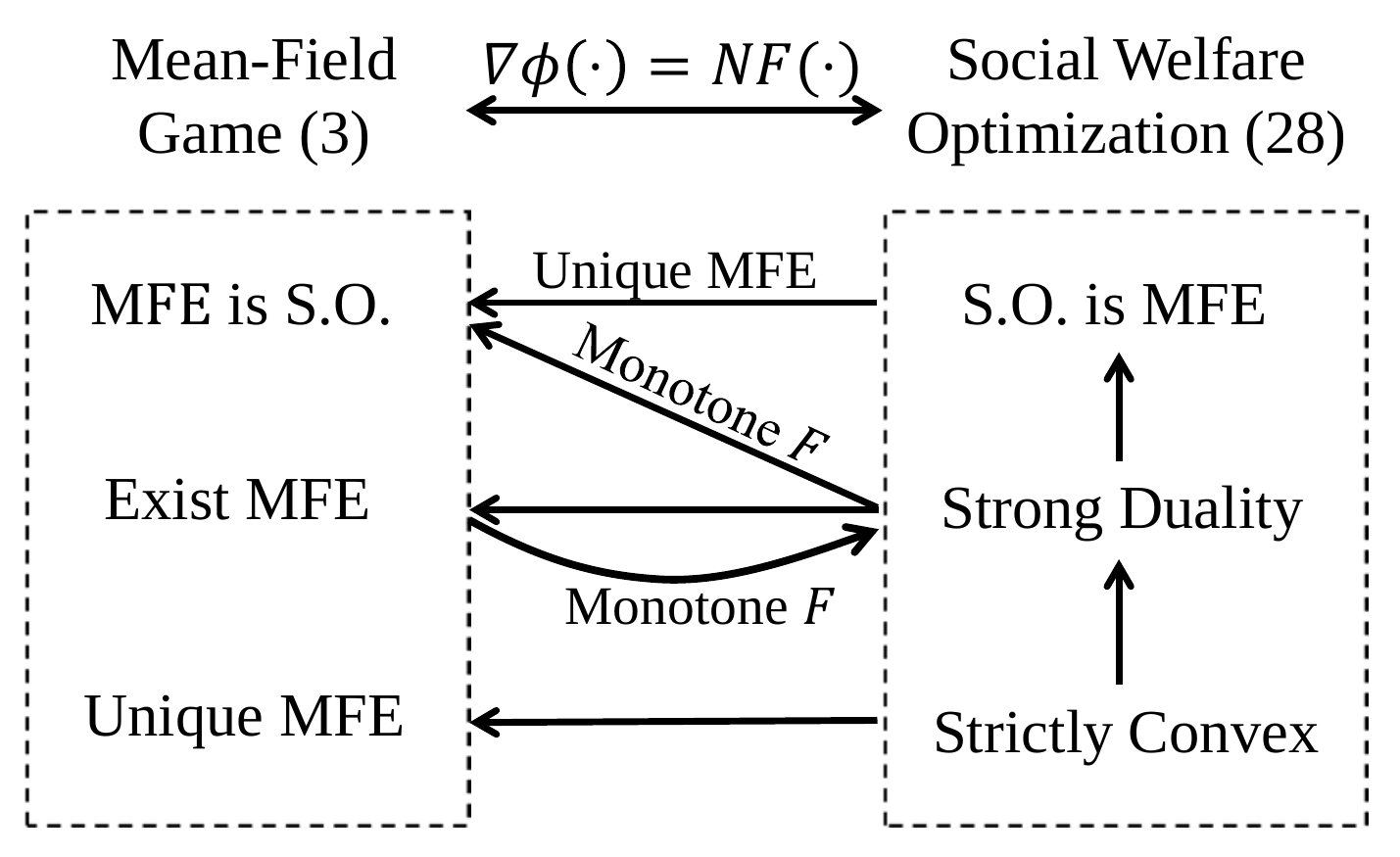}
\caption{The pictorial summary of our results, where S.O. stands for socially optimal solutions, and  MFE stands for mean-field equilibrium.}
\label{diagram}
\end{figure}

\subsection{Computation of the Mean-Field Equilibrium}
Aside from existence and uniqueness, another important implication of our result is on the computation of the mean-field equilibrium. Since the mean-field game can be connected to the social welfare optimization problem (\ref{socialwelfare_constructed}), we can compute the mean-field equilibrium by solving the corresponding social welfare optimization problem. When problem (\ref{socialwelfare_constructed}) is convex, there are many efficient algorithms to compute its solutions. In this subsection, we will present a primal-dual algorithm \cite[Chap. 10]{luenberger1997optimization} to compute the mean-field equilibrium, and we will show that many algorithms in the literature for computing the mean-field equilibrium are equivalent to the primal-dual algorithm for solving the corresponding social welfare optimization problem.

We consider a mean-field game (\ref{individualoptimization}) that satisfies Assumption \ref{convexityassumption}.  According to Corollary \ref{corollary1}, to compute the mean-field equilibrium, we can construct the corresponding social welfare optimization problem and solve it using a primal-dual algorithm. The details of the algorithm is summarized in Algorithm \ref{primaldualalgorithm}. To implement the algorithm, we first construct the social welfare optimization problem (\ref{socialwelfare_constructed}) by finding the virtual cost $\phi$. After this, an initial guess for the Lagrangian multiplier $\lambda$ is broadcast to all the agents.  Each agent can then independently solve the convex optimization problem (\ref{primaldual1}), while the virtual supplier solves the optimization problem (\ref{primaldual2}) for a given $\lambda$. The solutions of the cost minimization problems are collected and used to update the dual $\lambda$ according to (\ref{primaldual3}). The updated dual variable is then broadcast to the agents again and this procedure is iterated until it converges. Based on \cite[Chap. 10]{luenberger1997optimization} and Corollary \ref{corollary1}, it is easy to prove that the algorithm converges to the mean-field equilibrium of (\ref{individualoptimization}).
\begin{proposition}
Algorithm \ref{primaldualalgorithm} converges to the mean-field equilibrium of (\ref{individualoptimization}) if  Assumption \ref{convexityassumption} is satisfied 
and the step-size $\iota^k$ satisfies $\lim_{k\rightarrow \infty}\iota_k=0$ and $\lim_{k\rightarrow \infty}\sum_{m=1}^k \iota _m=\infty$.
\end{proposition}

\begin{algorithm}[!hbt]
\caption{The Primal-Dual Algorithm to Compute the Mean-Field Equilibrium} \label{primaldualalgorithm}
\begin{algorithmic}[1]

\REQUIRE the mean-field game (\ref{individualoptimization}),

\STATE  Construct (\ref{socialwelfare_constructed}) by finding $\phi(\cdot)$ that $\nabla \phi(\cdot)=NF(\cdot)$,

\STATE Generate initial guess for the Lagrange Multiple,  $\lambda^0$,

\FOR {$k=1,2, \dots,$}

\STATE Update the individual decisions by solving:
\begin{equation}{\hspace{-0.5cm}}
\label{primaldual1}
u_i^{k}= \argmin_{u_i\in \bar{\mathcal{U}}_i} \mathbb{E} \left(V_i(f_i(u_i,\pi_i),u_i)+\lambda^{k-1} \cdot f_i(u_i,\pi_i) \right)     
\end{equation}

\STATE Update the virtual supplier decision:
\begin{equation}
\label{primaldual2}
z^{k-1} =\argmin_{z\in \mathcal{X}} \phi(z)-N\lambda^{k-1} \cdot z ,     
\end{equation}

\STATE Update the dual variable according to:
\begin{equation}
\label{primaldual3} 
\lambda^{k} =\lambda^{k-1}+\iota^k  \left( \dfrac{1}{N}\mathbb{E} f_i(u_i^k,\pi_i)-z^{k-1} \right),
\end{equation}

\ENDFOR

\ENSURE the collective decisions $(u_1,\ldots,u_N)$.
\end{algorithmic}
\end{algorithm}


It can be verified that the proposed algorithm includes many existing ways to compute the mean-field equilibrium as special cases. For instance, the algorithms proposed in  \cite{rammatico2016decentralized} and  \cite{ma2013decentralized} are equivalent to the primal-dual algorithm with a scaled stepsize in (\ref{primaldual3}). Specifically, a finite-horizon deterministic linear quadratic mean-field game was considered in \cite{rammatico2016decentralized}, and an iterative algorithm was proposed to compute its mean-field equilibrium. The first step of the algorithm solves the same optimal response problem as (\ref{primaldual1}), and the second step updates $z$ according to:
\begin{equation}
\label{meanupdateinliterature}
z^{k}=z^{k-1}+\nu^k\left(  \dfrac{1}{N} \sum_{i=1}^N f_i(u_i^k) - z^{k-1} \right),
\end{equation}  
where $\nu^k$ is the step size. We comment that this is a scaled version of the  Algorithm \ref{primaldualalgorithm}. This is because in the linear quadratic game, $\phi(\cdot)$ is a quadratic function, and (\ref{primaldual2}) indicates that there is a positive definite matrix $B$ such that $Bz^k=\lambda^k$ for all $k$. Therefore, the update equation (\ref{primaldual3}) in Algorithm~\ref{primaldualalgorithm} can be written as:
\begin{equation*} 
z^k =z^{k-1}+B^{-1} \iota^k  \left( \dfrac{1}{N}\mathbb{E} f_i(u_i^k,\pi_i)-z^{k-1} \right),
\end{equation*}
which is equivalent to (\ref{meanupdateinliterature}) with $B^{-1} \iota^k=\nu^k$.  

\begin{figure*}[bt]%
\begin{minipage}[b]{0.32\linewidth}
\centering
\includegraphics[width = 0.93\linewidth]{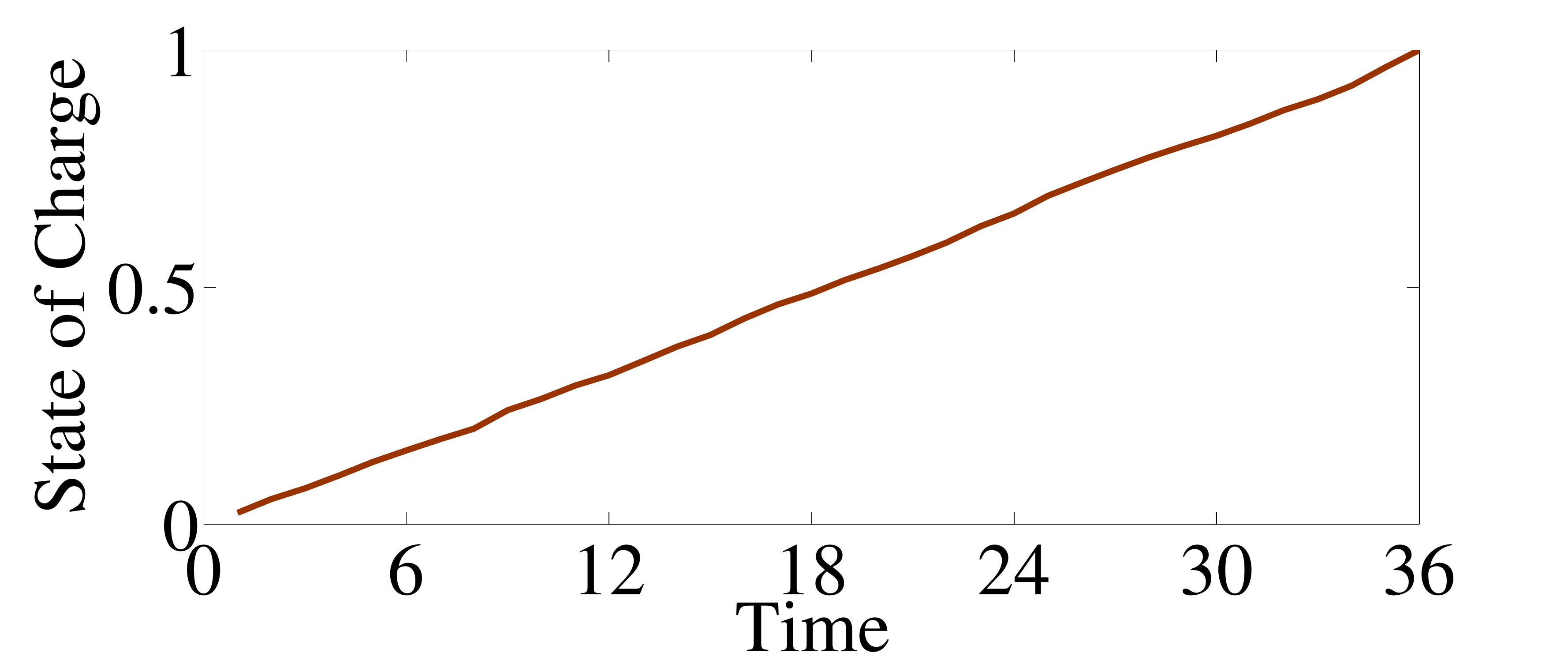}
\caption{The SOC of a randomly selected electric vehicle under ADMM.}
\label{fig1}
\end{minipage}
\begin{minipage}[b]{0.01\linewidth}
\hfill
\end{minipage}
\begin{minipage}[b]{0.32\linewidth}
\centering
\includegraphics[width = 0.95\linewidth]{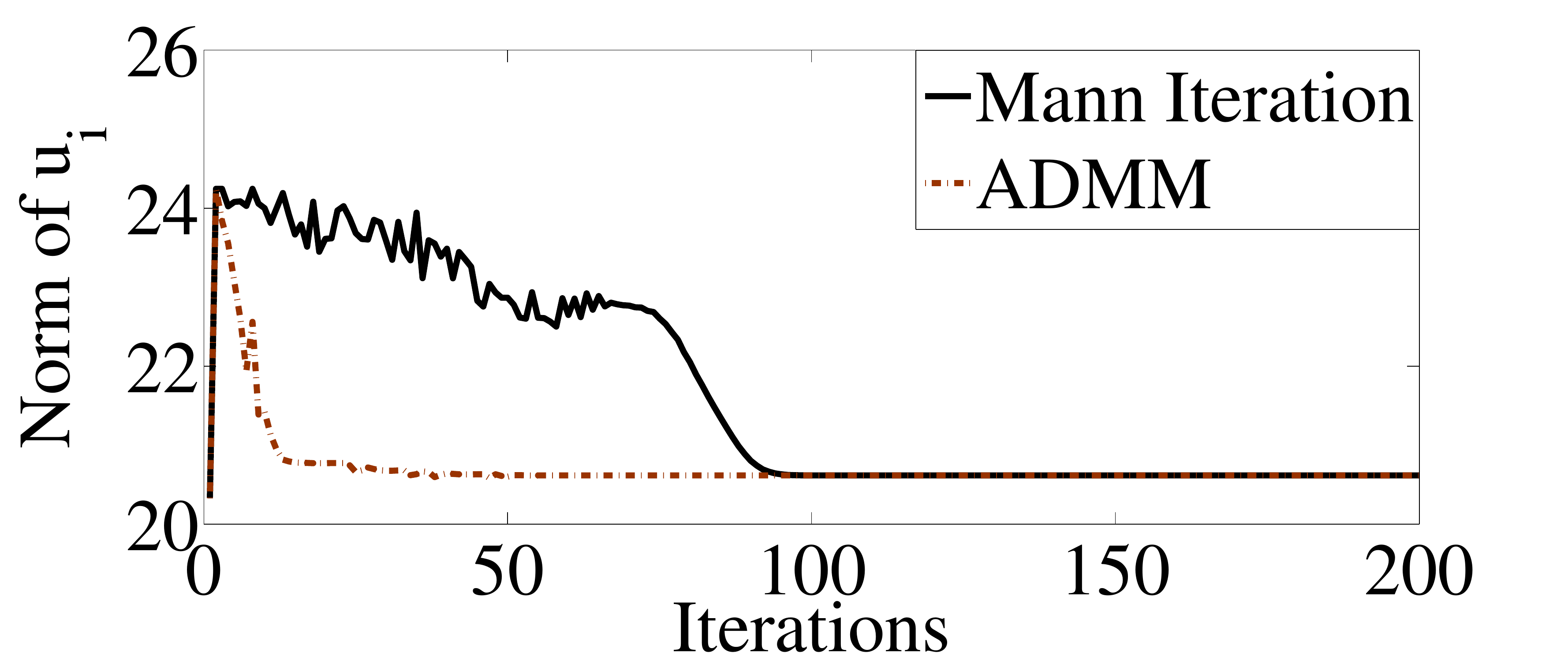}
\caption{The control decision of a randomly selected electric vehicle.}
\label{fig2}
\end{minipage}
\begin{minipage}[b]{0.01\linewidth}
\hfill
\end{minipage}
\begin{minipage}[b]{0.32\linewidth}
\centering
\includegraphics[width = 0.95\linewidth]{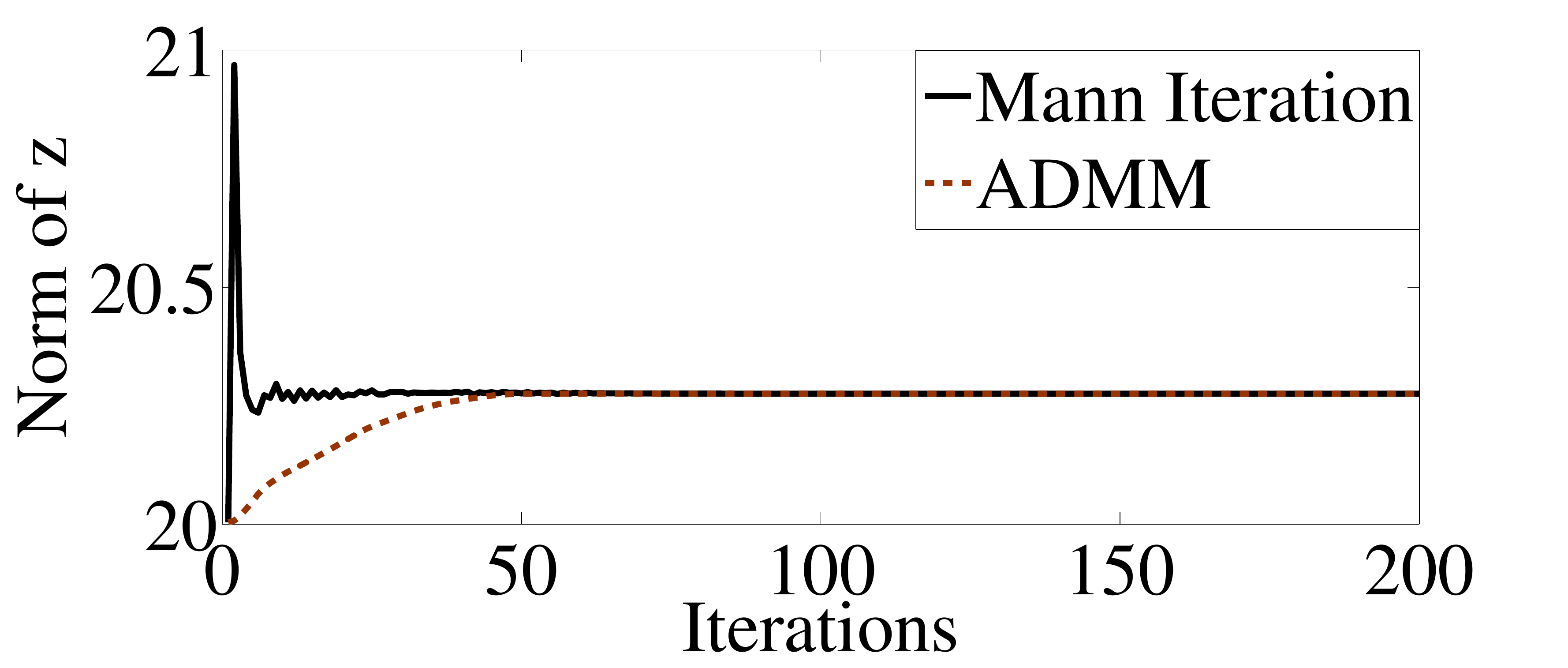}
\caption{The average control decision of the population in two algorithms.}
\label{fig3}
\end{minipage}
\end{figure*}


\section{Case Studies}
This section presents two examples to show how our results can be used to study the properties of the mean-field equilibrium. The first example is a special case where the mean-field coupling term $F(\cdot)$ is monotone. In this case, we can show that the mean-field equilibrium is equivalent to socially optimal solution if the social welfare optimization problem has strong duality.  The second example presents a more general case where the monotonicity condition is not satisfied. In this case, we can still draw connections between the mean-field equilibrium and the social welfare optimization problem using the result of Theorem \ref{mainresult1} and Corollary \ref{corollary2}.

\subsection{Example: Special Case with Monotone Mean-Field Coupling}
The first example considers the problem of coordinating the charging of a population of electric vehicles (EV) \cite{rammatico2016decentralized}. Each EV is modeled as a linear dynamic system, and the objective is to acquire a charge amount within a finite horizon while minimizing the charging cost. The charging cost of each EV is coupled through the electricity price, which is an affine function of the average of charging energy. This leads to the following game problem \cite{rammatico2016decentralized}:
\begin{align}
\label{PEVchargingame}
&\min_{u_i}  \eta||u_i-z||^2+2\gamma(z+c)^T u_i  \\
&\text{ subject to: } \nonumber\\
&\begin{cases}
   z=\dfrac{1}{N}\sum_{i=1}^N u_i, \nonumber\\
  x_i(t+1)=x_i(t)+u_i(t), \nonumber\\
  0\leq x_i(t) \leq \bar{x}_i, 0\leq u_i(t) \leq \bar{u}_i, \sum_{t=1}^Tu_i(t)=\gamma_i, \nonumber
\end{cases}
\end{align}
where $0< \eta \ll \gamma$, $x_i(t)\in \mathbb{R}$ is the state of charge (scaled by the capacity) of the EV battery, $u_i(t)$ denotes the charging energy during the $t$th control period, $\bar{x}_i$ is the battery capacity, and the electricity price is $2\gamma(z+c)$. 

In this mean-field game, the coupling term $F(z)=2(\gamma-\eta)z$ is increasing with respect to $z$, thus the monotonicity condition is satisfied. We also note that (\ref{PEVchargingame}) is slightly different from (\ref{individualoptimization}) in the sense that the coupling term depends on the average of control instead of the state. However, based on Remark \ref{remarkoncontrolmean}, there is no essential difference between these two formulations, and our result applies universally.

The focus of this case study is to investigate the efficiency, existence, uniqueness, and computation of the mean-field equilibrium to (\ref{PEVchargingame}). To this end,  we first construct the social welfare optimization problem for (\ref{PEVchargingame}) by finding the function $\phi(\cdot)$ such that $F(z)=\dfrac{1}{N} \nabla \phi(z)$. Since $F(z)=2(\gamma-\eta)z$, we have $\phi(z)=N(\gamma-\eta)z^Tz$. This indicates that the social welfare optimization problem for (\ref{PEVchargingame}) is as follows:
\begin{align}
\label{simulation_social}
&\min_{(u_1,\ldots,u_N,z)} \sum_{i=1}^N \left(\eta||u_i||^2+2\gamma c^Tu_i\right)+N(\gamma-\eta)z^Tz  \\
&\text{ subject to: } \nonumber\\
&\begin{cases}
z=\dfrac{1}{N}\sum_{i=1}^N u_i \\
x_i(t+1)=x_i(t)+u_i(t), \quad \forall i\in \mathcal{I},\quad \forall k\in \mathcal{K} \\
0\leq x_i(t) \leq \bar{x}_i, 0\leq u_i(t) \leq \bar{u}_i, \sum_{t=1}^Tu_i(t)=\gamma_i, \nonumber.
\end{cases}
\end{align}
where $\mathcal{I}=\{1,\ldots,N\}$ and $\mathcal{K}=\{1,\ldots,K\}$. It can be verified that the game problem (\ref{PEVchargingame}) satisfies Assumption \ref{convexityassumption}. According to Corollary \ref{corollary1}, the control decision $(u_1,\ldots,u_N)$ is a mean-field equilibrium for  (\ref{PEVchargingame}) if and only if it is an optimal solution to 
(\ref{simulation_social}). In addition, since (\ref{simulation_social}) is strictly convex with respect to $(u_1,\ldots,u_N, z)$, based on Theorem~\ref{uniquenesstheorem}}, the mean-field equilibrium exists and is unique.

Next, we study the computation of the mean-field equilibrium to (\ref{PEVchargingame}). In \cite{rammatico2016decentralized}, the mean-field equilibrium of (\ref{PEVchargingame}) is computed by an iterative algorithm.  In the algorithm, each agent takes $z$ as given and solve the problem (\ref{PEVchargingame}) to obtain an optimal control $u_i$. Based on $u_i$, $z$ is updated according to (\ref{meanupdateinliterature}). Each agent then takes $z'$ as given to repeat the first step.  This algorithm converges to the mean-field equilibrium of (\ref{PEVchargingame}) if $\lim_{k\rightarrow \infty}\nu_k=0$ and $\lim_{k\rightarrow \infty}\sum_{m=1}^k \nu_m=\infty$. Under this choice of step-size, the algorithm is referred to as the Mann iteration, and we will use it as the benchmark

In this paper, we propose to compute the mean-field equilibrium using the connections between the mean-field game and the social welfare optimization problem. 
Based on Corollary \ref{corollary1}, the mean-field equilibrium can be computed by solving the convex problem  (\ref{simulation_social}).
This problem can be efficiently solved by alternating direction method of multipliers (ADMM) \cite{boyd2011distributed}. In the rest of this subsection, we will use numerical simulation to compare the performance of the proposed ADMM algorithm with the benchmark algorithm (Mann iteration). 

In the numerical simulation, we generate 100 sets of EV parameters over 36 control periods, and each period spans 5 minutes. The heterogeneous parameters, including the capacity of the EV batteries, the maximum charging rate of the batteries, and the other parameters in the objective function are all generated based on uniform distributions. We run both Mann iterations and ADMM for 200 iterations, and the simulation results are shown in Fig.~\ref{fig1} - Fig. \ref{fig3}. In Fig. \ref{fig1}, we randomly select an EV and show its trajectory of state of charge under the ADMM solution. To compare the performance between ADMM and the Mann iteration, we show $||u_i||$ for a randomly selected EV (Fig. \ref{fig2}) and the average control decision $||z||$ (Fig. \ref{fig3}) over each iteration. Based on the simulation results, the ADMM algorithm convergences to the optimal solution after about 50 iterations, while the Mann iteration converges after 100 iterations. Therefore,  the proposed ADMM converges faster than Mann iteration. We emphasize that the algorithm converges only if both $z$ and $u_i$ converge. In Fig. \ref{fig3}, it seems that the Mann iteration converges faster than ADMM. As a matter of fact, although $||z||$ converges fairly fast, $z$ does not converge until after 100 iterations. This has been verified in the codes.

\subsection{Example: Non-Monotone Mean-Field Coupling}
The second example presents a general case where the mean-field coupling term is not monotone. In the rest of this subsection, we  show how to use the result of Theorem \ref{mainresult1} and Corollary \ref{corollary2} to draw connections between the mean-field equilibrium and the social welfare optimization problem. 

Consider a game with $N$ agents. Each agent $i$ wants to minimize the following objective function:
\begin{equation}
\label{examplemfg}
\min_{x_i} ||x_i||^2 +\kappa \int_0^{\infty} e^{-(\rho+1)t}x_i(t) \mbox{sin} \big(\bar{x}(t)\big)  dt,
\end{equation}
where $x_i$ takes value in the Hilbert space $\mathcal{X}=\{ x | \int_{[0,\infty)}e^{-\rho t}|x(t)|^2dt<\infty\}$ with the inner product (\ref{innerproduct4example}), $\bar{x}(t)=\dfrac{1}{N}\sum_{i=1}^N x_i(t)$,  and $\rho$ and $\kappa$ are positive scalars. In this example, the system dynamics are $x_i=u_i$, and  the mean-field coupling term is $F(\bar{x})=\{F_t(\bar{x}), t\geq 0\}$, where $F_t(\bar{x})=\kappa e^{-t}\mbox{sin}\big(\bar{x}(t)\big)$.
It is clear that the monotonicity condition does not hold.

Next, we will construct the corresponding social welfare optimization problem for (\ref{examplemfg}). To this end, consider a virtual supply cost that satisfies $\nabla \phi(z)=NF(z)$:
\begin{equation}
\label{virtualcost4example}
\phi(z)=-N\kappa \int_0^{\infty} e^{-t}\mbox{cos}\big(z(t)\big) dt,
\end{equation}
This gives rise to the following social welfare optimization problem:
\begin{align}
\label{socialwelfare_example}
&\min_{x_1,\ldots,x_N,z} \sum_{i=1}^N ||x_i||^2-N\kappa \int_0^{\infty} e^{-t}\mbox{cos}\big(z(t)\big)dt\\
&\text{ s.t.}
\begin{cases}
z=\dfrac{1}{N}  \sum_{i=1}^N x_i; \\
x_i=\mathcal{X}, z\in \mathcal{X},  \quad \forall i=1,\ldots,N.
\end{cases}
\end{align}
The above problem is non-convex with respect to $(x_1,\ldots,x_N, z)$, but we can still show that the primal problem  social (\ref{socialwelfare_example}) has the same optimal value as its dual problem. 
\begin{lemma}
\label{strongduality4example}
The social welfare optimization problem (\ref{socialwelfare_example}) has strong duality.
\end{lemma}
\vspace{-0.5cm}
\begin{pf*}{Proof.}
To prove strong duality, according to Definition \ref{strongdualitydefition}, it suffices to show that there exists $\lambda^*$ such that $P^*=D^*=D(\lambda^*)$. 
To show this, we first note that since $\kappa>0$, the cost function of primal problem (\ref{socialwelfare_example}) is lower bounded by $-N\kappa$:
\begin{equation*}
\sum_{i=1}^N ||x_i||^2-N\kappa \int_0^{\infty}e^{-t} \mbox{cos}\big(z(t)\big)dt \geq -N\kappa.
\end{equation*}
It can be verified that when $x_i=0$ and $z=0$, the cost function of (\ref{socialwelfare_example}) equals $-N\kappa$ and the constraints are satisfied. Therefore, $-N\kappa$ is the optimal value for the primal problem. 
According to Definition \ref{strongdualitydefition}, it suffices to find $\lambda^*$ such that the minimum value of $L(u,z,\lambda^*)$ is also $-N\kappa$. Let $\lambda^*=0$, then the Lagrangian dual $L(u,z,0)$ corresponds to the following problem:
\begin{equation}  
\label{dualproblem4example}
\min_{x_1,\ldots,x_N,z} \sum_{i=1}^N ||x_i||^2-N\kappa \int_0^{\infty}e^{-t} \mbox{cos}\big(z(t)\big)dt
\end{equation}
\quad \quad  s.t.: $x_i=\mathbb{R}, z\in \mathbb{R},  \quad \forall i=1,\ldots,N.$ \\
The optimal value of (\ref{dualproblem4example}) is clearly $-N\kappa$. Therefore, strong duality holds.
\end{pf*}

Based on Theorem {\ref{mainresult1}, if the social welfare optimization problem (\ref{socialwelfare_example}) has strong duality, then any solution to (\ref{socialwelfare_example}) is a mean-field equilibrium. 

We note that this relation only holds from one direction: there may exist a mean-field equilibrium which is not socially optimal. In section IV, we showed that this can be resolved when the mean-field equilibrium is unique. In our example, we can show that this is true under some conditions:
\begin{lemma}
\label{uniqueness4example}
If $0<\kappa<2$, then the mean-field equations for (\ref{examplemfg}) admit a unique solution.
\end{lemma}
\vspace{-0.5cm}
\begin{pf*}{Proof.}
To prove the uniqueness, the idea is to construct a contraction mapping, whose fixed point is the solution to the mean-field equations. In particular, we first regard $z$ as given and solve the problem (\ref{examplemfg}) to derive the optimal solution as $x_i^*(t)=-\dfrac{1}{2}\kappa \mbox{sin}\big(\bar{x}(t)\big)$. Then we define a function $T:\mathcal{X}\rightarrow \mathcal{X}$ that maps $\bar{x}$ to the average of $x_i^*$: $T(\bar{x})=\{T_t(\bar{x}),t\geq 0\}$, and $T_t(\bar{x})=-\dfrac{1}{2}\kappa \mbox{sin}\big(\bar{x}(t)\big)$.
It can be verified that the mean-field equilibrium is the fixed point of this mapping. 

Since $|\mbox{sin}(x)-\mbox{sin}(y)|\leq |x-y|$ for $x, y\in \mathbb{R}$, we have $||T(m_1)-T(m_2)||\leq \dfrac{1}{2}\kappa ||m_1-m_2||$ for any $m_1\in \mathcal{X}$ and $m_2\in \mathcal{X}$. Therefore, as $\kappa<2$, $T(\bar{x})$ is a contraction mapping, and it has a unique fixed point. This completes the proof.
\end{pf*}
Based on Corollary \ref{corollary2} and Lemma \ref{uniqueness4example}, it is clear that the mean-field equilibrium to  (\ref{examplemfg}) is socially optimal when $0<\kappa<2$.

To summarize, the connection between the mean-field equilibrium to (\ref{examplemfg}) and the social welfare optimization problem is as follows:
\begin{proposition}
If $\kappa>0$, the any socially optimal solution of (\ref{socialwelfare_example}) is a mean-field equilibrium to (\ref{examplemfg}). In addition, if $0<\kappa<2$, then the mean-field equilibrium to (\ref{examplemfg}) is the globally optimal solution to the social welfare optimization problem (\ref{socialwelfare_example}).
\end{proposition}
The proof follows from Lemma \ref{strongduality4example}, Lemma \ref{uniqueness4example}, Theorem \ref{mainresult1} and Corollary \ref{corollary2}.

\section{conclusion}
This paper studies the connections between a class of mean-field games and the social welfare optimization problem. We showed that the mean-field equilibrium is the optimal solution to a social welfare optimization problem, and this holds for both convex and non-convex individual cost functions and action spaces. The result enables us to evaluate and promote the efficiency of the mean-field equilibria, and it also provides interesting implication on the existence, uniqueness and computation of the mean-field equilibrium. Numerical simulations are presented to validate the proposed approach. Future work includes extending the proposed approach to the case of infinitely many agents and more general formulations where the mean-field term depends on the probability distribution of the population state.

\section*{Acknowledgement}
The authors would like to thank Prof. Tamer Basar (UIUC) for fruitful discussions, and Prof. Roland Malhame (PolyMtl) and Prof. Sergio Grammatico (TU Delft) for their insightful comments and helpful suggestions.

\section*{Appendix}

\subsection*{\bf{A: Proof of Lemma \ref{hilbertspace}}}
It is well-known that based on the Riesz-Fischer theorem \cite[p.148]{royden1988real}, the $L^2$ space is complete. We show that $\mathcal{X}$ is isomorphic to $L^2$, and therefore it is also complete \cite[p.20]{conway2013course}. For this purpose, we define the following mapping: $g:L^2\rightarrow \mathcal{X}$ that satisfies $g(\cdot)=\{g_t(\cdot),t\geq 0\}$ and $g_t(l(t))=e^{\rho t/2}l(t)$ for any $l\in L^2$, where $l=\{l(t),t\geq 0\}$. This is a linear surjective mapping and it can be verified that $g(l_1(t))\cdot g(l_2(t))=l_1(t)\cdot l_2(t)$, where the left-hand side inner product is defined as  (\ref{innerproduct4example}) on $\mathcal{X}$ and the right-hand side inner product is defined on the $L^2$ space in the canonical form. This indicates that $L^2$ and $\mathcal{X}$ are isomorphic, which completes the proof. 

\vspace{0.1cm}

\subsection*{\bf{ B: Proof of Lemma \ref{lemma1}}}
To prove this result, it is clear that we have:
\begin{align*}
\big|\mathbb{E} \big(F(\bar{x})\cdot x_i\big)-&F(\mathbb{E}\bar{x})\cdot \mathbb{E} x_i \big|=\big| I_1+I_2 \big| \leq \nonumber\\ 
&\big| I_1\big| +\big| I_2\big|,
\end{align*}
where we define $I_1$ as $I_1=\mathbb{E} \big(F(\bar{x})\cdot x_i\big)-\mathbb{E} F(\bar{x})\cdot\mathbb{E} x_i$ and $I_2=\mathbb{E} F(\bar{x})\cdot\mathbb{E} x_i-F(\mathbb{E}\bar{x})\cdot \mathbb{E} x_i$. Then it suffices to show that both $\big| I_1\big|$ and $\big| I_2\big|$ converge to $0$ at the rate of $\dfrac{1}{\sqrt{N}}$. To show this, we first note that:
\begin{align*}
\big| I_1\big|=\big| \mathbb{E} \big(F(\bar{x})\cdot x_i\big)-&\mathbb{E} F(\bar{x})\cdot\mathbb{E} x_i \big|=\big| I_3+I_4\big| \leq \nonumber \\
&\big| I_3\big| +\big| I_4\big|,
\end{align*}
where we define $I_3=\mathbb{E} \big(F(\bar{x})\cdot x_i\big)-\mathbb{E}\big(F(\bar{x}_{-i})\cdot x_i\big)$, $I_4=\mathbb{E} \big(F(\bar{x}_{-i})\cdot x_i\big)-\mathbb{E}\big(F(\bar{x})\cdot \mathbb{E}x_i\big)$ and $\bar{x}_{-i}=\dfrac{1}{N}\sum_{j\neq i}x_j$.
Since $F(\cdot)$ is Lipschitz continuous with the constant $L\geq 0$, and the second moment of $x_i$ is bounded, we have:
\begin{align*}
\big| I_3\big|=\big| \mathbb{E}& \big(F(\bar{x})\cdot x_i\big)-\mathbb{E}\big(F(\bar{x}_{-i})\cdot x_i\big) \big| \leq \nonumber\\
\mathbb{E}\big| F(\bar{x})\cdot x_i&-F(\bar{x}_{-i})\cdot x_i\big| \leq \mathbb{E} \big(\big\Vert \dfrac{L}{N} x_i\big\Vert \big\Vert x_i  \big\Vert \big) = \nonumber \\
&\dfrac{L}{N}\mathbb{E}||x_i||^2 \leq \dfrac{LC}{N}. 
\end{align*}
In addition, as $x_i$ is uncorrelated with $\bar{x}_{-i}$, we have 
\begin{align*}
\big| I_4\big|=\big| \mathbb{E} \big(F(\bar{x}_{-i})\cdot x_i\big)-\mathbb{E}\big(F(\bar{x})\cdot \mathbb{E}x_i\big) \big| = \nonumber\\
\big|\mathbb{E} F(\bar{x}_{-i})\cdot \mathbb{E}x_i-\mathbb{E}F(\bar{x})\cdot \mathbb{E}x_i \big| \leq \dfrac{L}{N}\big\Vert \mathbb{E}x_i\big\Vert^2. 
\end{align*}
Note that $\mathbb{E}||x_i||^2$ is bounded, and thus $\big\Vert\mathbb{E}x_i\big\Vert$ is also bounded:
\begin{align}
\label{boundonexpectation}
\big\Vert\mathbb{E}x_i\big\Vert \leq \mathbb{E} \big\Vert x_i\big\Vert &= \sqrt{ \mathbb{E} \big\Vert x_i\big\Vert ^2- \mathbb{E} \left( \big\Vert x_i\big\Vert-\mathbb{E} \big\Vert x_i\big\Vert  \right)^2 } \nonumber \\
&\leq \sqrt{\mathbb{E} \big\Vert x_i\big\Vert ^2 } \leq \sqrt{C}.
\end{align}
This indicates that $\big| I_4\big|\leq \dfrac{LC}{N}$.
 Therefore, $\big| I_1\big|\leq \big| I_3\big|+\big| I_4\big|\leq \dfrac{2LC}{N}$. To show that $\big| I_2\big|$ converges to $0$ at the rate of $\dfrac{1}{\sqrt{N}}$, we define a random variable $r_N=F(\bar{x})-F(\mathbb{E}\bar{x})$. Note that:
 \begin{equation}
\label{tempresultinproof1}
I_2\leq \big\Vert \mathbb{E}F(\bar{x})-F(\mathbb{E}\bar{x})\big\Vert\big\Vert \mathbb{E}x_i \big\Vert = \big\Vert\mathbb{E} r_N\big\Vert \big\Vert \mathbb{E}x_i \big\Vert.
\end{equation} 
Since $||\mathbb{E} r_N||\leq \mathbb{E}||r_N||$ and $||\mathbb{E}x_i||\leq \sqrt{C}$,  the inequality (\ref{tempresultinproof1}) reduces to $I_2\leq \sqrt{C}\mathbb{E}||r_N||$.
Therefore, for our purpose, it suffices to show that $\mathbb{E}||r_N||=O(\dfrac{1}{\sqrt{N}})$.
Since $F(\cdot)$ is Lipschitz continuous, we have:
\begin{equation*}
||r_N||=||F(\bar{x})-F(\mathbb{E}\bar{x})||\leq L||\bar{x}-\mathbb{E}\bar{x}||.
\end{equation*}
Therefore, it suffices to show that $\mathbb{E}||\bar{x}-\mathbb{E}\bar{x}|| =O(\dfrac{1}{\sqrt{N}})$. 
To prove this, we note that since $x_i$ and $x_j$ are uncorrelated, thus we have the following relation:
\begin{align*}
\mathbb{E}||\bar{x}-\mathbb{E}\bar{x}||^2&=\dfrac{1}{N^2}\mathbb{E}\sum_{i=1}^N ||x_i-\mathbb{E}x_i||^2< \\
&\dfrac{1}{N^2}\mathbb{E}\sum_{i=1}^N ||x_i||^2\leq \dfrac{NC}{N^2}=\dfrac{C}{N}.
\end{align*}
Therefore, similar to (\ref{boundonexpectation}), we have 
\begin{equation*}
\mathbb{E}||\bar{x}-\mathbb{E}\bar{x}||\leq \sqrt{\mathbb{E}||\bar{x}-\mathbb{E}\bar{x}||^2} \leq \sqrt{\dfrac{C}{N}}.
\end{equation*}
This completes the proof.
\vspace{0.1cm}

\subsection*{\bf{ C: Proof of Lemma \ref{lemma2} }}
Since $\nabla G(\cdot)$ is Lipschitz continuous with constant $\beta$, for any $x,y\in \mathcal{X}$, we have:
\begin{equation}
\label{lipschitzgradient}
|G(x)-G(y)|\leq \nabla G(y) \cdot (x-y)+\dfrac{\beta}{2}||x-y||^2.
\end{equation}
Based on (\ref{lipschitzgradient}), the following relation holds:
\begin{align}
\label{betasmooth2}
\mathbb{E} |G(\bar{x})-G(\bar{x}_{-i})|&\leq \mathbb{E} \left( \nabla G(\bar{x})  \cdot \dfrac{-x_i}{N} \right)+\dfrac{\beta}{2}\mathbb{E} ||\dfrac{x_i}{N}||^2 \nonumber \\
&\leq \mathbb{E} \left( ||\nabla G(\bar{x})||||\dfrac{x_i}{N}|| \right)+   \dfrac{\beta}{2N^2} \mathbb{E} ||x_i||^2.
\end{align}
Since derivative of $G(\cdot)$ is Lipschitz continuous with constant $\beta$, we have $||\nabla G(\bar{x})||\leq ||\nabla G(0)||+\beta||\bar{x}||$. Plugging this in (\ref{betasmooth2}), we obtain:
\begin{align*}
\mathbb{E} |G(\bar{x})-G(\bar{x}_{-i})|\leq & 
\dfrac{||\nabla G(0)||}{N}\mathbb{E}||x_i||+
\dfrac{\beta}{N} \mathbb{E} \left( ||\bar{x}|| ||x_i|| \right) \\
&+   \dfrac{\beta}{2N^2} \mathbb{E} ||x_i||^2
\end{align*}
The right-hand side of the above inequality consists of three terms. We will show that all the three terms converge to 0 at the rate of $1/N$. For notation convenience, let the three terms be $T_1^N=\dfrac{||\nabla G(0)||}{N}\mathbb{E}||x_i||$, $T_2^N=\dfrac{\beta}{N} \mathbb{E} \left( ||\bar{x}|| ||x_i|| \right)$ and $T_3^N=\dfrac{\beta}{2N^2} \mathbb{E} ||x_i||^2$. In the proof of Lemma~\ref{lemma1}, we showed that that $\mathbb{E}||x_i||\leq \sqrt{C}$. Therefore, we have $T_1^N<\dfrac{\sqrt{C}||\nabla G(0)||}{N}$. Due to Assumption \ref{assumption_cost}, $||\nabla G(0)||$ is bounded, thus $T_1^N=O(\dfrac{1}{N})$. As for the second term, we have:
\begin{align*}
T_2^N=\dfrac{\beta}{N} \mathbb{E} \left( ||\bar{x}|| ||x_i|| \right)= \dfrac{\beta}{N^2} \mathbb{E} \left( \sum_{j=1}^N||x_j|| ||x_i|| \right) \\
=\dfrac{\beta}{N^2} \left( \mathbb{E} ||x_i||^2  +    \sum_{j\neq i} \mathbb{E} (||x_i||||x_j||)   \right) \\
=\dfrac{\beta}{N^2} \left( \mathbb{E} ||x_i||^2  +    \sum_{j\neq i} \mathbb{E}||x_i|| \mathbb{E}||x_j||   \right) \leq  \dfrac{\beta C}{N}
\end{align*}  
Furthermore, since $\mathbb{E}||x_i||^2\leq C$, we have $T_3^N=\dfrac{\beta}{2N^2} \mathbb{E} ||x_i||^2\leq \dfrac{\beta C}{2N^2}$, thus $T_1^N+T_2^N+T_3^N=O(\dfrac{1}{N})$. This completes the proof.

\vspace{0.1cm}

\bibliographystyle{unsrt}
\bibliography{TransactiveControl}

\end{document}